# Pierre van Hiele, David Tall and Hans Freudenthal: Getting the facts right

Thomas Colignatus
http://thomascool.eu
November 3 & 4 2015

**Abstract**

Pierre van Hiele (1909-2010) suggested, both in 1957 and later repeatedly, wide application for the Van Hiele levels in insight, both for more disciplines and for different subjects in mathematics. David Tall (2013) suggests that Van Hiele only saw application to geometry. Tall claims that only he himself now extends to wider application. Getting the facts right, it can be observed that Tall misread Van Hiele (2002). It remains important that Tall supports the wide application of Van Hiele's theory. Tall apparently didn't know that Freudenthal claimed it too. There appears to exist a general lack of understanding of the Van Hiele - Freudenthal combination since 1957. Hans Freudenthal (1905-1990) also suggested that Van Hiele only saw application to geometry, and that only he, Freudenthal, saw the general application. Freudenthal adopted various notions from Van Hiele, misrepresented those, gave those new names of himself, and started referring to this instead of to Van Hiele. The misrepresentation may clarify why Tall didn't recognise Van Hiele's theory. Freudenthal mistook Van Hiele's distinction of concrete versus abstract for the distinction of reality versus model (applied mathematics). Freudenthal's misconception of "realistic mathematics education" (RME) partly doesn't work and the part that works was mostly taken from Van Hiele. This common lack of understanding of Van Hiele partly explains the situation in the education in mathematics and the research on this. Another factor is that mathematicians like Freudenthal and Tall are trained for abstraction and have less understanding of the empirics of mathematics education.

## *Introduction*

Pierre van Hiele (1909-2010) & Dina (Dieke) Van Hiele-Geldof (1911-1958) were inspired by Jean Piaget's idea of levels of understanding of mathematics, notably linked to age. They tested that idea, and empirically developed and defined the Van Hiele levels of insight, more independent of age. I will also refer to them as levels of abstraction in understanding mathematics. There are two separate theses from July 5 1957 – the year of Sputnik. Pierre stated that the idea of the levels was his. As Piaget presented a general theory, the Van Hieles presented an alternative general theory. The usual reference is to Van Hiele (1986) in English. I will refer to Van Hiele (1957 & 1973). Pierre's statement that insight wasn't necessarily related to age for example had an impact on later Russian education on mathematics.

Pierre van Hiele had been teaching mathematics on various topics for two decades and had observed the levels in various topics, and indeed in various disciplines, like chemistry and didactics itself. His thesis merely took geometry as an example, or, rather as the example *par excellence*, as geometry is a foundation stone for mathematics as the art and science of *demonstration*. His thesis provided in mathematical fashion both a *definition* of the levels in human understanding and an *existence proof*. See **Appendix A** also on the proper translation of the title of the thesis.

David Tall (2013), *"How humans learn to think mathematically. Exploring the three worlds of mathematics",* suggests that Pierre van Hiele had a limited understanding of the portent of this theory, notably that Van Hiele saw it limited to geometry and not applicable for e.g. algebra or even other disciplines. Tall also suggests that he himself extends our understanding to the wider portent of those levels. However, it appears that Tall's suggestions are based upon a misreading of Van Hiele (2002). The truth is that Van Hiele was quite aware of the fundamental nature of his and Dieke's result. It is better to get the facts right and indeed alert students of education (not only mathematics) to the wealth that can be found in Van Hiele's work. In itself it is important that Tall supports the notion that the Van Hiele levels have general validity.

Tall apparently wasn't aware of the similar claim by Freudenthal either. Hans Freudenthal (1905-1990) was Pierre's thesis supervisor – with second advisor M.J. Langeveld (1905-1989). It appears that Freudenthal also claimed at a later moment that Van Hiele saw only limitation to geometry, and that it was he, Freudenthal, who found the general applicability. The topic is important enough for a full chapter in Sacha la Bastide – Van Gemert (2006) – henceforth LB-VG – which is a Dutch thesis on Freudenthal. This thesis inconsistently first allows Van Hiele's claim and then proceeds showing that it was Freudenthal who at a later moment discovered this general applicability.

Updating this 2014 article gave the option to write about Tall and Freudenthal separately.

(1) However, the corroboration by Tall (2013) of the theory is important. It will be useful for readers to know that we are dealing with a serious issue. We will regard mathematics as dealing with abstraction. The Van Hiele theory concerns the general didactics of abstraction and abstract concepts. Abstraction is often seen as 'higher' and possibly more complex, but it is better seen as the elimination of aspects in concrete elements, so that it may also be regarded as rather simple. In this way we can also understand that the human brain is capable of abstract thought, as leaving out aspects doesn't have to be much of an achievement. The Van Hiele levels of insight are as important for epistemology as the law of conservation of energy is for physics, see Colignatus (2015). (It is an open question, for example, whether Harvard's Jeanne Chall's (1921-1999) stages in development of reading are related to the Van Hiele levels.)

(2) Below discussion will not recall what the Van Hiele theory of levels actually is. This is a deliberate choice. One of the objectives of this paper is to alert the reader to the wealth of insight in didactics in the work by the Van Hieles, as applicable to didactics in general, and thus it is consistent to refer only. There are more expositions available in English, see Van Hiele (1986) or **Appendix A** on Van Hiele (1959) at ERIC. Regrettably, Pierre's thesis never got translated into English. The discussion of Tall (2013) will cause a partical restatement of the Van Hiele levels anyway. Thus, by discussing both Tall and Freudenthal in one paper, we have the advantage of common exposition, terminology and references. And we may better understand how Freudenthal's terminology might have caused Tall (and others) not to recognise Van Hiele's theory in disguise.

(3) The update of this 2014 article in 2015 concerns only that new sources have become available for readers of English. The facts and conclusions on Freudenthal are already in the 2014 version. Readers of English can now benefit from the following:

(a) There is a 2015 English translation of LB-VG. This allows for an independent check on my own translations in 2014 from that Dutch 2006 thesis in **Appendix B**.
(b) All editions of Euclides since 1924 have come available on the internet. This journal of the Dutch association of teachers of mathematics allows international readers to verify Van Hiele's 1957 claim. Perhaps the Van Hieles theses didn't get wide circulation, and Van Hiele (1986) in English is apparently out of print, but their 1957 article in Euclides now is available to the world. A section below will give the translation of the claim. Also, Freudenthal (1948) is an early discussion on didactics for highschool (with a contribution by Van Hiele in the same volume). Colignatus (2015b) shows that Freudenthal (1948) is rather traditional and fits the period. Translation can be tricky. When Freudenthal (1948) uses the phrases *"far advanced stage of abstraction"* and *"higher point of view"* then this may remind of the Van Hiele levels. However, he only refers to a vantage point, like one can climb a hill to see further. The evidence lies in that Freudenthal doesn't discuss level transitions. Van Hiele really arrived at a new contribution, turning Piaget's basic notion of stages into logically defined categories that are relevant for both theory and practice.
(c) I provided a 2015 English translation myself on a section of the Alberts & Kaenders (2005) interview with Van Hiele (then 96 years of age), Colignatus (2015g).

A question by a reader was why the critique w.r.t. Freudenthal cannot be restricted to saying that he should have referred more. For this reason there now is a separate section on Freudenthal's breach in research integrity that focuses on four points of evidence. One can verify that this information is already in the 2014 version of this article.

**Putting Van Hiele into the geometry box**

There is indeed some tendency, like on wikipedia 2014, to restrict the Van Hiele levels to geometry alone but it wouldn't be right and useful if this became a general misunderstanding. By analogy, when someone mentions the number 4 as an example of an even number, it is invalid to infer that this person thinks that 4 gives all even numbers.

The impact of misconception can be large. Apparently the levels are not applied to other subjects much. The task group on learning processes of the 2008 US National Mathematics Advisory Panel states, see Geary et al. (2008:4:xxi-xxii):

> "The van Hiele model (1986) has been the dominant theory of geometric reasoning in mathematics education for the past several decades. (...) Research shows that the van Hiele theory provides a generally valid description of the development of students' geometric reasoning, yet this area of research is still in its infancy."

Professor Tall's personal reappraisal of Van Hiele's work is important for the wider recognition of that work. The convincing part in Tall's argumentation for the wider application is basically no different from the argumentation that Van Hiele already provided, which is another corroboration of the original insight. A nuance is that research in the education of mathematics has provided additional corroboration since 1957 also by Van Hiele himself. Another nuance is that Tall (2013) unfortunately still has a limited understanding of the Van Hiele levels and introduces some misunderstandings. While this personal reappraisal is important, there might be the danger that researchers and students now would focus on Tall (2013) as the most recent text while it would be advisable to study the original work by the Van Hieles.

Updating this article caused me to read a section again from Freudenthal (1987), the autobiography at age 82 (online). On p354-355, my emphasis, we see again that he puts Van Hiele in the geometry box while Freudenthal as thesis supervisor should know that Pierre presented a general theory:

> "Like I took the **Van Hiele's geometric levels** and interpreted them as applying to mathematics as a whole, the activity that takes place at a lower level becomes subject of research at the higher level. At the higher level the mathematical activity at the lower level becomes object of relection: by reflecting about your mathematical activity you create new mathematics – which is done by the original discoverer and the learning re-discoverer." [1]

Presmeg (2014:52-53) recalls – and confirms in an email that Van Hiele's theory is wider than only geometry, and that the original research was in the 1950s (Colignatus (2015i)):

---

[1] "Zoals ik Van Hieles meetkundige niveaus op de wiskunde als geheel toegepast interpreteerde, wordt wat op het lager niveau wiskundig handelen was, onderwerp van onderzoek op het hogere - op het hogere wordt de wiskundige activiteit van het lager niveau object van reflectie: door op je wiskundig handelen te reflecteren schep je nieuwe wiskunde - de eerste ontdekker én de lerende herontdekker."

"One vivid memory I have, which is pertinent to this topic, concerns the interview process when I applied for the professorship that was open in the Mathematics Department at ISU in 2000. During the two days of interviews, one of the meetings was with the mathematicians. They asked me what I, as a mathematics education researcher, could offer them in their work. I thought quickly, and then described the theoretical model of levels of learning geometry put forward by the van Hieles as a result of their research in The Netherlands in the 1970s. The mathematicians could see the value of such research in teaching mathematics. I was hired!"

**Van Hiele and Tall in 2014, joined by Freudenthal, and new evidence in 2015**

This article has two layers. Sections with (2014) are copied from Colignatus (2014d) that originally focused on Tall's misconception. At a few points I insert 2015 comments in square brackets. New sections use the news in 2015. Colignatus (2014d) has been on the internet for a year and has had a role in some discussions. Leaving these sections intact allows readers to understand part of that discussion in 2014-2015. **Appendix B** could be moved to the main body of the text, but because of the first layer it is left in its 2014 position. The 2015 English translation of LB-VG should allow readers to verify that this 2014 **Appendix B** was correct.

Tall (2013) is the most recent publication and it is from a world-renowned researcher who after retirement takes stock of his lifetime work. He is also at some distance from the Van Hiele and Freudenthal interaction in Holland. It is useful to start with Tall (2013) and his supposed evidence in Van Hiele (2002), and give the quotes that highlight both the claim and the misunderstanding. Subsequently, we provide quotes from Dutch sources now in English translation. Subsequently we compare the Van Hiele levels with Tall's diagram of such levels, to check that we are speaking about the same things. Some of Tall's misunderstandings generate a somewhat distorted model, which causes difficult semantics whether we are really speaking about the same things: but overall this sameness could be accepted when the misunderstandings are recognized for what they are.

Subsequently, we look at Freudenthal, with the 2015 English access to the evidence. The enlightening observation in 2014 already was that also Freudenthal claimed a general theory derived from Van Hiele, and also provided distorted information as well. This seems to have affected Tall's perception as well. **Appendix B** contains supporting translations from a Dutch thesis by La Bastide-Van Gemert (2006) of which Chapter 7 looked at Freudenthal on the Van Hiele levels.

Let us however start with the basic evidence from 1957 and Van Hiele's own view in the interview by Alberts & Kaenders (2005) (with my 2015 translation).

## *Basic evidence from 1957*

The basic evidence is in the thesis by Pierre van Hiele (1957). This however is not generally available. The Van Hieles (1957) Euclides article provides this statement on p45, that you can retrieve online now in 2015 and e.g. submit to Google Translate. My own translation and emphasis is, see Colignatus (2015c):

> "Above, we presented a didactic approach to introducing geometry. This approach has the advantage that students experience how you can make a field of knowledge accessible for objective consideration. For such a field a requirement is that students already have command of global structures. They experience how they proceed from those to further analysis. **The approach presented here for geometry namely can be used also for other fields of knowledge (disciplines).** Whether it will be possible to treat such a field also in mathematical manner depends upon the nature of the field. For mathematical treatment it is necessary, amongst others, that the relations do not lose their nature when they are transformed into logical relations. For students, who have participated once in this approach, it will be easier to recognize the limitations than for those students, who have been forced to accept the logical-deductive system as a ready-made given. Thus we are dealing here with a formative value (Bildung), that can be acquired by the education in the introduction into geometry." (p45) [2]

In this translation, I adopt the English preference for shorter sentences, and I insert the word "Bildung" to better express what the Van Hieles intend. Naturally I am advising to have an independent translation, by an expert who is aware of the pitfalls in terminology in mathematics education research (like the Freudenthal *vantage point*). (We also see Stellan Ohlsson's view that students start from global / vague ideas (which isn't higher abstraction), see Colignatus (2015ad).)

## *The 2005 interview*

In 2014 I asked interviewers Alberts & Kaenders (2005) whether they could provide for an independent English translation of their interview with Van Hiele in 2005. Since they didn't do so, Colignatus (2015g) gives my translation of some key parts, reproduced here. Naturally I still call for an independent translation of the whole interview. For readers of English it is new information in 2015 that this part of the interview has a translation now.

It is useful for the reader to be aware that the 2014 version of this article was written with my knowledge of the interview in Dutch. I don't know whether someone provided David Tall with an independent translation. For some readers it came as a shock to see that Freudenthal and Van Hiele were no good friends.

---

[2] Dutch original with abbreviations replaced, and "h" included in "mathematiseren": "De hiervoor aangeduide wijze om het meetkunde-onderwijs te beginnen heeft het voordeel, dat de leerlingen ervaren, hoe men een kennisgebied, waarvan men globale structuren bezit, door analyse voor objektieve beschouwingen toegankelijk kan maken. De hier voor de meetkunde aangegeven weg kan namelijk ook voor andere kennisgebieden gebruikt worden. Of het daar ook mogelijk zal zijn het kennisveld tenslotte te mathematiseren, hangt van de aard van het veld af. Noodzakelijk daarvoor is immers onder andere, dat de relaties niet gedenatureerd worden, wanneer zij in logische relaties worden omgezet. Voor hen, die aan deze werkwijze eens aktief hebben deelgenomen, zal het gemakkelijker zijn de grenzen te herkennen dan voor hen, die het logisch deduktieve systeem als een kant en klaar gegeven hebben moeten aanvaarden. We hebben hier dus te doen met een vormende waarde, die verkregen kan worden door het onderwijs in het begin van de meetkunde."

The risk of the interview is that some readers may start regarding Van Hiele (VH) as an old man with a grudge, perhaps with a bad memory of real events (96 in 2005), who accuses Freudenthal (F) of misdeeds who cannot defend himself. Shouldn't VH and F have discussed this around say 1980 and have settled this as gentlemen ? Perhaps F would have been shocked to hear the criticism by VH ? But VH says that F was a bossy person. In other interviews about F as teacher and collegae, see Verhoef & Verhulst (2010), we see it confirmed that F was bossy indeed, also towards ministers of education. Even today researchers find it difficult to do something about it when their work is maltreated by others. We can indeed establish objectively that F abused work by VH: and it is easy to check, see below, once you are alerted to it. The interview of 2005 shows as a whole that VH was clear of mind. Holland has more experience with people of 96 years of age who can recall issues accurately (and other issues not, of course). VH doesn't express a grudge but makes a matter of fact observation. It is really a pity that there is no translation of the full interview, so that one can check that this protest is only a part of a longer and interesting discussion about mathematics education.

Copying from Colignatus (2015g): In 2005, Gerard Alberts (mathematician, historian) and Rainer Kaenders (mathematician, educator) interviewed Pierre van Hiele (1909-2010). The interview was published in the journal of the Royal Dutch Society for Mathematics, as G. Alberts & R. Kaenders (2005), *"Interview Pierre van Hiele. Ik liet de kinderen wél iets leren"*, NAW 5/6 nr. 3, september, p247-251. The publication is in Dutch and I will translate some parts into English.

The introduction to the interview is:

> "Pierre van Hiele is the silent force in didactics of mathematics in The Netherlands. He was teacher of mathematics and chemistry and never much stepped in the floodlights. His work receives broad international recognition and one cannot think about didactics of mathematics without it. His work is still studied, amongst others *Stucture and Insight*. Van Hiele is ninety-six." [3]

On page 247:

> "My relation with Freudenthal wasn't so good, that I would go and drink coffee with him. Besides, Freudenthal has later frequently sabotaged my work, guys." [4]

Page 251 on Hans Freudental:

> "What role did Freudenthal play in your life? "I did not mix well with Freudenthal. From the beginning. He was a bossy person. He did cause me to get ideas. That is rather all."

---

[3] Dutch: "Pierre van Hiele is de stille kracht van de didactiek van de wiskunde in Nederland. Hij was wiskunde- en scheikundeleraar en is nooit veel op de voorgrond getreden. Zijn werk vindt brede internationale erkenning en is tegenwoordig niet meer weg te denken uit de wiskundedidactiek. Zijn werk, waaronder het invloedrijke *Begrip en inzicht, werkboek van de wiskundedidactiek* wordt nog steeds bestudeerd. Van Hiele is zesennegentig jaar."

[4] "Zo goed was mijn relatie met Freudenthal niet dat ik met hem ging koffie drinken. Trouwens, Freudenthal heeft mij later nogal eens een hak gezet, jongens."

Freudenthal used different descriptions of the process of abstraction. In the *Vorrede zu einer Wissenschaft vom Mathematikunterricht* he presented this process in terms of comprehension and apprehension. Did he also think differently about the role of levels of insight? "Yes, I believe actually that he did not really understand much about the levels of insight."

Freudenthal was your thesis supervisor (promotor). Did he also help you in stepping outside of the small circle – with contacts outside of The Netherlands? "The last thing definitely not. No, the situation was actually that I had to vouch for myself. For example I remember a conference in America, at which a speaker referred to my work and said: 'Mr. Van Hiele whom I am mentioning now is actually present in this very lecture hall. Mr. Van Hiele, please rise (so that everyone can recognize you)." Someone in the audience, a German, asked where he could read about my work. I replied that there would appear a book of mine in English shortly. Then Freudenthal who was also present said: 'You can also read about it in my book.' Which wasn't true. He just was sabotaging me again. Freudenthal was like this, yes." He was sabotaging you all the time? "Actually yes. Freudenthal never was a friendly person for me, no." Where can he have been sabotaging you? Didn't you work in entirely different environments? "Yes, but he tried to pinch something from me all the time."

Later you got more recognition for your work. Were you able to make peace with him then? "Well, peace? No, actually not. In that case you first would have made war. I don't make war."

You presume that Freudenthal did not fully understand your work. Did you understand him, conversely? "Yes, I understood what I knew of him. And I often agreed with it too. I wasn't in constant quarrel with Freudenthal. From his side, he had very much respect for my ideas on vectors in primary education. He praised me very much for that.""  [5]

## *Breach of research integrity*

Let us distinguish the courtroom from the world of ideas. The facts presented here may not be evidence enough to cause a verdict according to some set of laws of one country or other. However, for the research in mathematics education it is important to infer from these facts that Freudenthal has been in breach of research integrity and created a false line of research that misdirected some generations of researchers. He adopted various notions from Van Hiele, misrepresented those, gave those new names of himself, and started referring to this instead of to Van Hiele: which is appropriation of work and the withholding of proper reference. The misrepresentation doesn't make the theft less so.

Various aspects have been mentioned over 2014-2015 and one may lose track. It will be useful to collect four main facts. These were already known in Colignatus (2014d), the former version of this paper. The news is only that sources become available for readers of English.

---

[5] The original can be found in Alberts & Kaenders (2005) and Colignatus (2015g), both online. For vectors, Colignatus (2015e) shows that pupils at elementary school can prove the Pythagorean Theorem using the method by Yvonne Killian. Colignatus (2015f) shows how set theory can be brought into the highschool programme without getting lost on Cantor's transfinites.

(1) While the Van Hieles focused on the process from concrete to abstract, Freudenthal substituted this with a process from reality to model (applied mathematics), forgetting or denying that one must first master mathematics before one can apply it. Freudenthal's misconception of "realistic mathematics education" (RME) partly doesn't work and the part that works was mostly taken from Van Hiele. Only an abstract thinker or non-didact might think to argue that Freudenthal's approach could pass for good didactics. Freudenthal's confusion is discussed in *"Conquest of the Plane",* Colignatus (2011a), Chapter 15.

New in 2015: (a) Since I am translating quotes, I also took some quotes from the 1987 autobiography, see **Appendix D**. (b) Psychologist Stellan Ohlsson argues that learning goes "from abstract to concrete" but when we remain with epistemological terminology then he means to say that it goes from vague to precise, see Colignatus (2015ad).

(2) LB-VG (2006:201) has this quote from Freudenthal (already in English), who acknowledges that the Van Hieles already **used** guided re-invention, **purposively,** so that Freudenthal by implication only provides phraseology:

> "It is not by chance that the Van Hieles seized upon this idea. To my knowledge they were the first who wrote a textbook in which the learning process is purposively initiated and kept up as a process of re-inventing. [ftnt 79]"

(3) Let me quote from **Appendix B**, quoting from LB-VG (2006):

> p182 gives a quote by Freudenthal in his autobiographic book p354, which is rather convoluted and lacks the clarity that one would expect from a mathematician:
>
> "The process of mathematisation that the Van Hieles were mostly involved with, was that of geometry, more exactly put: they were the first who interpreted the geometric learning process as a process of mathematisation (even though they did not use that term, [6] and neither the term re-invention). In this manner Pierre discovered in the educational process, as Dieke described it, the levels of which I spoke earlier. I picked up that discovery - not unlikely the most important element in my own learning process of mathematics education." [ftnt] [7]
>
> Comment: Freudenthal thus suggests: (a) Pierre's insight is just seeing what Dieke described, so that she would be the real discoverer. (b) Freudenthal's words "mathematisation" and "re-invention" would be crucial to describe what happens in math education, otherwise you do not understand what math education is about, and it is only Freudenthal who provided this insight. (c) The Van Hieles wrote about geometry but were limited to this, so that it was Freudenthal himself who picked it up and provided the wider portent by means of his new words.

---

[6] Above 1957 quote from Euclides shows that the Van Hieles did use the term "mathematise" (though without a h), meaning "to treat a subject in mathematical manner".
[7] Freudenthal (1987:354), (at age 82),
http://www.dbnl.org/tekst/freu002schr01_01/freu002schr01_01_0025.php

A reaction to this by readers has been that this analysis is not necessary, and that Freudenthal's quote can also be read differently. One might argue that Freudenthal with "mostly" doesn't exclude that the Van Hieles were doing other things. This reaction doesn't account for the fact that Freudenthal was a mathematician. He would prefer to be accurate – unless there was a (subconscious) factor that caused some hiding.

Let me give the reasoning of 2014 in smaller steps. This seems to be the crux of the argument:

> *Can we expect mathematicians to be accurate in their statements ? With a proper distinction between a general theory and a particular example ?*

Presmeg referred to the Van Hiele levels in geometry but did not intend to review the work by the Van Hieles. Instead, someone who wants to give an **accurate summary** of the Van Hiele 1957 thesis and theory of levels of insight:

- **does not state**: "they were the first who interpreted the geometric learning process as a process of mathematisation"
- **instead states**: "they were the first who interpreted the general process of learning and teaching mathematics as a process of mathematisation, demonstrating this by the introduction into geometry"

If Freudenthal was a mathematician then the inaccuracy cannot have been other than deliberate. Thus: *he deliberately gave an inaccurate summary*. This amounts to a deliberate misrepresentation. There is also an appropriation, given his claim of his "own" learning process (and publications about this). Inescapably: this is fraud.

Again from **Appendix B**:

> p194, taking a quote from Freudenthal's autobiographic book p352:
>
> "What I learned from the Van Hieles I have reworked in my own manner - that is how things happen." [8]
>
> Comment: This is the veiled confession of appropriation. Freudenthal claims to be powerless and innocent of deliberate appropriation since "that is how things happen". Who however considers what that "reworking" involves sees only phraseology and lack of proper reference.

(4) There is Van Hiele's protest, quoted from the Alberts & Kaenders (2005) interview. Van Hiele denies that his theory is presented well in Freudenthal's books. Still, Freudenthal claimed this in public discussion. One might argue that Freudenthal would have believed so honestly, and that he believed honestly that he improved on the approach. The above however shows that there was appropriation and fraud, so that this "honest belief" was self-delusory. The Van Hiele interview is proof for the systematic and public abuse. Van Hiele deserves some credit as witness of what has happened.

---

[8] My translation of: "Wat ik van de Van Hieles leerde heb ik op mijn eigen wijze verwerkt – zo gaat dat nu eenmaal."

**On the issue of intent**

As said: this article is not directed at the courtroom but is concerned with research and the world of ideas. Still, there are some readers who seem to suggest that the analysis in this paper would be vague when it doesn't *prove* whether Freudenthal was merely confused or had deliberate intent. This is a non-sequitur.

One would suppose that there are hard-nosed detectives who will argue that it is a profession of itself to recover intent. As a researcher I have had no training on this, and thus my findings are just what they are. I am no historian either, and only look at these matters from the angle of mathematics education research. I want to refer to the proper sources, and should be able to distinguish whether ideas differ or are merely given different names.

(1) Van Hiele stated the general applicability of his theory of levels. For mathematicians the distinction between a general idea and an example is obvious. That Freudenthal in later years put Van Hiele in the box of geometry only and claimed the discovery of the general applicability for himself, is *for a mathematician* deliberate intent at misrepresentation, and thus fraud.

(2) Giving new names to what the Van Hieles did, without checking with them whether this was okay and that he got it right, and stopping to refer to them, is appropriation. One might argue that this is what bossy persons do honestly, but it still is improper.

(3) These conclusions don't change when Freudenthal was confused on other aspects. Freudenthal didn't really understand the Van Hiele theory of levels. Van Hiele stated this, and it is indeed shown in the confusion of "concrete versus abstract" with "reality versus model" (applied mathematics). There is an argument that one cannot appropriate an idea when one doesn't understand it. This argument is false, since a crow may steal a wedding ring without knowing what it represents. Freudenthal didn't state that he couldn't understand the theory of levels, and must have thought that he did understand it. He also seems to have understood the general applicability. It seems likely, though we will never be sure, that Freudenthal had no intent to misunderstand Van Hiele: thus honestly misunderstood him. This still is misrepresentation but without intent. Freudenthal's association of education with applied mathematics is fitting for a mathematician who has been trained on abstraction. It gives the development of "realistic mathematics education" (RME) a distinct difference with Van Hiele didactics. It may suggest own creativity but actually is the exercise of claimed authority in a field that one isn't trained for. Nowadays this will be regarded as another form of fraud, though the meddling of mathematicians in the education of mathematics might be an eternal phenomenon.

These musings may have only historical value if we restrict the framework to these old men. The full framework however is mathematics education research. It helps to be aware that RME has false foundations and that Van Hiele provides useful foundations. In Holland, psychologist Ben Wilbrink wants to discard Van Hiele's theory because he regards it as part and parcel of RME that doesn't work, see Colignatus (2015a), but one should keep proper perspective.

## *Quoting Tall and Van Hiele (2014)*

Tall (2013:153): "As we consider the whole framework of development of mathematical thinking, we see a prescient meaning in the title of van Hiele's book *Structure and Insight.* [footnote referring to Van Hiele (1986)] Even though he saw his theoretical development of levels of structure applying only to geometry and not algebra [footnote referring to Van Hiele (2002)], his broad development, interpreted as structural abstraction through *recognition, description, definition* and *deduction,* can now be extended to apply throught the whole of mathematics." (Note the limitation to "structural" abstraction. This comes back below.)

Tall (2013:430): "It was only in 2011, when Pierre van Hiele passed away at the grand old age of 100, that I explicitly realized something that I had 'known' all along: that the structural abstraction through *recognition, description, definition* and *deduction* applied successively to the three worlds as concepts in geometry, arithmetic and algebra, and formal mathematics were recognized, described, defined and deduced using appropriate forms of proof."

When I queried professor Tall on this, he sent me a copy of Van Hiele (2002) so that it has been verified that he regards that article as the "proof" that Van Hiele would have only a limited understanding of the portent of his theory. It appears to be a misunderstanding, and I don't think that such "proof" could be found elsewhere either. We can constrast above two quotes from Tall (2013) with Van Hiele (2002).

Van Hiele (2002:46) in his conclusion, with the word "disciplines" referring to also physics, chemistry, biology, economics, medicine, language studies, and so on:

> "In most disciplines there are different levels of thinking: the visual level, the descriptive level and the theoretical level."

Indeed, Van Hiele gives various examples in my copy of *"Begrip en Inzicht"* (1973), which I presume will be more extended in English in *"Structure and Insight"* (1986) that however is not in my possession.

Thus, Van Hiele was aware of the wide portent of the theory of levels of abstraction.

How could it be, that Tall did come to think otherwise ? One aspect will be an issue of reading well. One aspect might be a general misconception that the theory applied only to geometry. Readers who suffer this misconception might no longer read carefully. It may be observed that geometry has somewhat been reduced in the education in mathematics so that, if there is a misconception that the Van Hiele levels apply only to geometry, then this is one avenue to explain the reduced attention for the Van Hiele levels of abstraction.

The following quotes are relevant for the view on algebra. Van Hiele (2002) does not give a technical development of particular levels for arithmetic and algebra, but the point is that it shows that Van Hiele thought about the teaching of those subject matters in terms of levels. It is a wrong reading by Tall not to recognise this.

This "part of algebra" should not be mistaken for all algebra. Van Hiele (2002:28) warns:

> "The problems in algebra that cause instrumental thinking have nothing to do with level elevation since the Van Hiele levels do not apply to that part of algebra. People applied terms such as 'abstraction' and 'reflection' to the stages leading from one level to the next. This resulted in a confusion of tongues: we were talking about completely different things."

Van Hiele (2002:39): "The transition from arithmetic to algebra can not be considered the transition to a new level. Letters can be used to indicate variables, but with variables children are acquainted already. Letters can be used to indicate an unknown quantity, but this too is not new." Indeed, Tall (2013:105) confirms: "It is well known that students have much greater success in solving an equation with $x$ only on one side. [footnote] Filloy and Rojano [footnote] named this phenomenon the 'didactic cut' between arithmetic and algebra." Thus note:

- Van Hiele uses "algebra" to mean that the use of a single variable would still be classified as "algebra" even though it isn't really different from arithmetic.
- Filloy and Rojana use "algebra" as distinctive from arithmetic, so that a single occurrence of an unknown would not be classified as "algebra". (Or they use a dual sense.)
- Apart from this issue in terminology, the diagnosis is the same.
- And this isn't evidence that Van Hiele had a limited view on the portent of levels, but rather the opposite.

Van Hiele (2002:43): "The examples Skemp mentions in his article about I2, R2 and L2 do not have any relations with a level transition. They are part of algebra in which topic, as I have emphasised before, normally level transitions do not occur." Again "part of algebra" should not be read as all algebra. In this case we must observe that the sentence can be bracketed in different ways: "They are (part of algebra) in which topic ..." or as "They are part of (algebra in which topic ....)". An observant reader will be aware of this issue and rely on the rest of the text to determine the proper bracketing. My diagnosis is that professor Tall focussed on this sentence and mislaid the "which". However, the other sentences and in particular the conclusion in the very same article on page 46 (Van Hiele (2002:46), quoted above) make clear that this shouldn't be done. I have asked professor Tall whether this particular misunderstanding was the source of his suggestions indeed, but haven't received an answer on this particular question yet. [Which is also the case when writing this in 2015.]

Thus, Van Hiele was aware of the portent of his theory, contrary to what Tall states.

One might hold that it doesn't matter whether Tall read something wrong and that the useful issue is to arrive at a working theory for didactics. However, the point is that Van Hiele already presented such a theory. The derived issue is quite limited here: that Tall creates a confusion and that it helps to eliminate the reasons upon which he based that confusion: (i) so that others do not follow that same road, (ii) so that readers of Tall (2013) are alerted to this confusion there. Overall, students of didactics are advised to consider the original Van Hiele texts.

## *Some other sources (2014)*

It will be useful to point to some other (Dutch) sources and publications by Van Hiele also explicitly on arithmetic and algebra.

The opening paragraph in the *Introduction* to his thesis Van Hiele (1957:vii):

> ""Insight" is a concept that can present itself to us in different fashions. The meaning of the different aspects differs, depending upon the context for which one studies insight. In the following study here I have occupied myself with the position, that insight takes in the context of didactics and even more special in the didactics of geometry. This limitation causes that the conclusions that are arrived at cannot be regarded without additional research as "generally valid". From this study it may however appear that there shouldn't be expected differences in principle between "insight in geometry" and "insight in mathematics in general". I am also under the impression that "insight in mathematics" will be quite similar in many respects with "insight in non-mathematical school subjects". It might however be that insight plays a much less fundamental role in some school subjects other than mathematics, without the implication that such school subjects would have to be less important for the child." [9]

Van Hiele (1957) is remarkably consistent in focussing on geometry, but discusses in passing that there are also levels of insight in algebra (p131-132) and in understanding didactics and insight itself (p201-204), while he also discusses Langeveld's study of learning checkers (p105). Overall, one cannot conclude that Van Hiele thought that his new theory of levels of insight was limited to geometry only.

La Bastide-Van Gemert (2006) - henceforth LB-VG - reports about Van Hiele on writing his thesis in 1957:

> "He restricted himself to the education in geometry, since he did not see differences in principle between insight in geometry and insight in mathematics in general." [10]

---

[9] My translation of: ""Inzicht" is een begrip, dat zich op verschillende wijzen aan ons kan voordoen. De betekenis van de verschillende aspekten varieert, al naar gelang men het inzicht in de ene of in de andere samenhang bestudeert. In de hier volgende studie heb ik mij speciaal bezig gehouden met de plaats, die het inzicht inneemt in de didaktische kontekst en nog meer speciaal in de didaktiek van de meetkunde. Deze beperking maakt, dat men de gevonden konkusies niet zonder nader onderzoek als "algemeen geldig" mag beschouwen. Uit deze studie moge echter blijken, dat er in ieder geval tussen "inzicht in meetkunde" en "inzicht in wiskunde in het algemeen" geen principiële verschillen verwacht mogen worden. Ook schijnt het mij toe, dat "inzicht in wiskunde" toch nog op vele punten overeenkomst zal vertonen met "inzicht in niet-wiskundige schoolvakken". Het zou echter wel kunnen zijn, dat in sommige schoolvakken het inzicht een minder fundamentele rol vervult dan in de wiskunde, zonder dat daardoor deze schoolvakken voor het kind minder belangrijk behoeven te zijn."

[10] My translation of: "Hij beperkte zich daarbij tot het meetkundeonderwijs, aangezien hij geen principiële verschillen zag tussen inzicht in meetkunde en inzicht in wiskunde in het algemeen." (p190)

She also reports that the two Van Hieles wrote in an article in *Euclides*, the Dutch journal for teachers of mathematics, in 1957 (they include some conditions that are not relevant here):

> "The approach given here for geometry can namely also be used for other disciplines." [11]

Van Hiele (1959) considers thought and then focuses on geometry as an example. His opening line is: "The art of teaching is a meeting of three elements: teacher, student, and subject matter." He speaks about mathematics in general rather than geometry.

Geometry enters only when: "The following example will illustrate what I mean." - with that example taken from geometry, all aware that the subject matter might affect the analysis. However he continues speaking about "optimal mathematical training" in general. A useful quote is: "In general, the teacher and the student speak a very different language. We can express this by saying: they think on different levels. Analysis of geometry indicates about five different levels."

It may be useful to know that Van Hiele's attention for the role of language derived from the teachings by Gerrit Mannoury (1867-1956), professor in mathematics at the university of Amsterdam in 1917-1937 who did much research in semiotics (or in Dutch: significa). Mannoury explained already quite early what Ludwig Wittgenstein rediscovered and rephrased more succinctly: *the meaning of a word is its usage* (i.e. the notion of *language games*). Van Hiele (1959) gives the closing statement: not about geometry only but for each discipline:

> "The heart of the idea of levels of thought lies in the statement that in each scientific discipline, it is possible to think and to reason at different levels, and that this reasoning calls for different languages. These languages sometimes use the same linguistic symbols, but these symbols do not have the same meaning in such a case, and are connected in a different way to other linguistic symbols. This situation is an obstacle to the exchange of views which goes on between teacher and student about the subject matter being taught. It can perhaps be considered the fundamental problem of didactics."

Van Hiele (1962) *"The relation between theory and problems in arithmetic and algebra"* is a chapter in an ICMI report by Freudenthal (ed) 1962. I haven't been able to check yet whether the levels are applied but would be amazed if they would not appear.

In their memorial text Broekman & Verhoef (2012:123) refer to Van Hiele (1964), a contribution in German to Odenbach (ed) 1964:

> "In the background there was the struggle by Van Hiele with working with the two different intuitions that already could be seen with Pythagoras and other ancients. When children learn, this can be seen in the intuition for continuity and that for discreteness, as this shows up in (spatial) geometry or counting (with integers) respectively. That discreteness concerns the transition from experience to abstractions in the form of symbols - thus to detach yourself from the image that is experienced and that determined the 'number'. In view comes arithmetic, numbers are nodes in a large network of relationships [reference]." [12]

---

[11] My translation of: "De hier voor de meetkunde aangegeven weg kan nl. ook voor andere kennisgebieden gebruikt worden." (p202)

All this should not surprise us. The Van Hieles started with Piaget's theory, which was a general theory of development. Their alternative was another general theory of development. They only took geometry as their test case for their theses.

W.r.t. the problem of induction: it is hard to prove a theory for all disciplines, even those not invented yet. The Van Hieles were aware of the limits of empirical methodology. See however Colignatus (2011c) for the *"definition & reality methodology":* definitions will be chosen such that reality can be covered with minimum uncertainty, and in some cases we might achieve virtually zero uncertainty. [See Colignatus (2015a) for a developed statement on this.] It is interesting to observe that Van Hiele (1957:191) was actually rather aware of this too:

> "We have mentioned already how the result of a study is often largely established by the choice of the definition." [13]

## *Diagram, embodiment and abstraction (2014)*

I do not intend to review Tall's book here, only to set the record straight w.r.t. Van Hiele. But perhaps a general remark is allowed. At points the reader can embrace Tall's objectives, yet at other points one wonders whether he has actually used Van Hiele's work in practice. Too often we see Tall perform as a mathematician trained for abstraction and too often we don't see Tall perform as an empirical scientist who has recovered from his training for abstraction. My diagnosis is that Tall (2013) is a seriously misdirected book and needs a full rewrite. I will explain this in more detail at another place. [Not yet done in 2015.] Let me now consider the issue of Tall's diagram of the "three worlds of mathematics", and the issues of embodiment and abstraction, as they relate to Van Hiele's levels of understanding of mathematics (i.e. abstraction).

Tall (2013:17-19) - see his online PDF of the first chapter - presents a diagram (or table) in which we can recognise the Van Hiele levels. Tall's new format introduces separate attention for the senses, notably vision and sound (language, symbolics). In principle also touch and motion would be important but this might be taken along in "language". One might test the aspects by using the vision and sound buttons on the tv-control, on broadcasts with or without subtitles. We essentially see the two hemispheres of the brain, with the prefrontal cortex monitoring. It might lead too far but it has been suggested that

---

[12] My translation of: "Op de achtergrond speelde hierbij de worsteling van Van Hiele met het werken met twee verschillende intuïties die ook al bij Pythagoras en zijn tijdgenoten te onderkennen waren. Dit komt voor lerende kinderen tot uiting in de intuïtie van continuïteit en die van discreetheid zoals die naar voren komen in de (ruimte-) meetkunde respectievelijk het (met gehele getallen) rekenen. Die discreetheid heeft betrekking op de overgang van aanschouwelijkheid naar abstracties in de vorm van symbolen — dus het zich losmaken van het aanschouwelijke beeld dat het 'aantal' bepaald heeft. Het rekenen komt in zicht, getallen zijn knooppunten geworden in een groot relatienet [reference]."

[13] My translation of: "Wij hebben er al eerder op gewezen, hoe door de keuze van de definitie het resultaat van het onderzoek dikwijls al grotendeels vastligt."

Greek culture was visual and Oriental culture was aural, so that we could understand Euclid's *Elements* as a result from the clash of civilizations in Alexandria: to write down in linguistic legal fashion what the visual mind could perceive. Overall, it seems a useful idea of Tall indeed to use the visual method of a diagram in two dimensions to display the field of discussion. Overall, this should be used with caution however, since symbolics like Roman XII or Hindu-Arabic 12 might be attributed to "language" but clearly have visual aspects. Problematic in Tall's diagram is his treatment of abstraction and his use of the term "embodied".

Above, we quoted Van Hiele (2002:28): "People applied terms such as 'abstraction' and 'reflection' to the stages leading from one level to the next. This resulted in a confusion of tongues: we were talking about completely different things." Van Hiele means that levels in understanding, with their web of relations of concepts, and actually diffferent meanings of the same words (actually speaking another language), cannot be merely reduced to such a vague term as "abstraction". When this is clear, I however would like to suggest that it can be advantageous to refer to the *levels of understanding of mathematics* (insight) as *levels in abstraction*. These are all somewhat vague notions while there is good reason to regard ["mathematics" as dealing with "abstraction"].

It is important to emphasize that thought is abstract by nature. When an apple creates an image in the mind (with all available senses, not just visual), then this mental image is abstract, and the brain can start processing it. For the education of mathematics it is crucial to be aware of this abstract nature of thought. It was the error of Hans Freudenthal to misunderstand the basic Van Hiele level: mistaking Van Hiele's reference to concreteness for some kind of "experience of reality", and even derive the name of his "realistich mathematics education" from it, and to introduce all kinds of real-world aspects into the curriculum and lose sight of that essential abstract nature of mathematics. It is the challenge for mathematics education to get pupils to focus on the abstract aspects that teachers know are useful to focus on. Colignatus (2011a) section 15.2 discusses how an essential step can be made here, in the "Conquest of the Plane". Colignatus (2011b) has some comments w.r.t. brain research with this role for abstraction.

Let us first present the Van Hiele levels of abstraction (insight) and then look at Tall's diagram. It is a sobering thought that we are basically classifying the math subjects of elementary school through university, but the thrust of the Van Hiele levels is the associated didactics, of providing the pupils with the appropriate materials and instructions. I feel a bit ashamed of presenting Table 1 but the confusion by professor Tall requires an answer, and the point remains that one better considers the wealth in the research by the Van Hieles for educational practice.

**Table 1**: **Van Hiele levels of abstraction (insight) distinguished by brain function**

| Deductive | 3. Formal | Euclid | | Hilbert |
|---|---|---|---|---|
| | 2. Informal | | School analytic geometry | School algebra |
| Practical | 1. Description | Elementary school Geometry | Overlap, such as procept | Arithmetic |
| | 0. Intuition | Visual emphasis | Overlap v & a | Aural emphasis |

In the basic Van Hiele level, the pupil is mainly sensing the world. Basically the outside world feeds the memory via the senses, as the memory itself isn't well developed to feed the mind. The pupil learns to recognise objects and the space in which the objects occur. My preference for the basic Van Hiele level is the term "intuition" as it indicates a key role for the subconscious mind as opposed to conscious cognition at the higher levels (even though the interaction between these two is more complex than passing on of sensory data). The point isn't quite this experience of reality, as Freudenthal might suggest water running from a faucet as the experience of a linear process. The point is that the experience is concrete, like a line drawn with a ruler, so that the process of abstraction has traction to start from something (close to the intended mental image). Van Hiele has the opposition concrete vs abstract, Freudenthal model vs reality, see Colignatus (2014c).

For level 3, the lawyers at the department of mathematics might argue that Euclid isn't quite formal enough to pass for the claim of being "formal" conform David Hilbert, since his definitions and axioms refer to notions in the "intended interpretation" of visual space. However, instead of including additional levels (as originally by Van Hiele), it remains a good idea of Tall to include these extra columns. Then, though Euclid contains a lot of words, he is conveniently put in the column with the visual emphasis.

Perhaps superfluously, it may be remarked that aspects of logic and set theory and notions of proof should already be taught at elementary school. It is the failure of mathematics educators with the Sputnik "New Math" and "realistic mathematics education" (RME), and so on, that causes that we do not develop what is potentially possible. In practice education at elementary school already uses notions of proof as it isn't all rote learning, so that above distinction between "deductive" and "practical" is again one of emphasis and degree and level of abstraction.

Let us now look at Tall's diagram. Tall (2013:17)'s "three worlds of mathematics" are: "One is based on (conceptual) embodiment, one on (operational) symbolism, and the third on (axiomatic) formalism, as each one grows from earlier experience." He also uses "structural, operational and formal abstraction". (Check the Tall (2013:153) quote above.)

**Diagram**: "Figure 1.5" copied from Tall (2013:17)

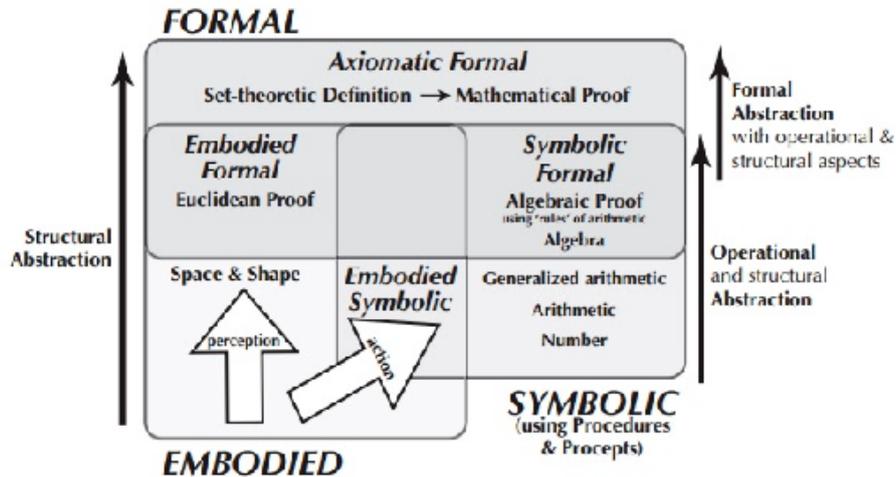

Figure 1.5: Preliminary outline of the development of the three worlds of mathematics

Professor Tall proposes to use the term "embodied". This term however has an inverse meaning w.r.t. the Van Hiele level. Van Hiele considers the situation that objects of the outside world create images in the mind of the pupil. An inverse process starts with an idea. We say that *an idea is embodied in something*. Tall must reason as a mathematician for whom a soccer ball that the pupil plays with is an embodiment of the mathematical idea of a sphere. Perhaps Tall's use of the word "embodied" is acceptable for objects, but it becomes awkward when the pupil learns the properties of empty space, (e.g. that a meter is the same in any direction); but we can stretch the meaning of embodiment too. In that case physical space embodies some math space. (Even though math space is *empty* by abstraction, and "embodiment" of *nothing* is a difficult notion to grasp, especially as some suggest that physical space isn't empty.)

Curiously though, Tall gives a definition of embodiment that focuses on the mental image, as in "the word became flesh", e.g. the pupil's mental concept of the soccer ball. Tall (2013:138-139) discusses President Bush senior's "decade of the brain" with various studies, and concludes: "The framework proposed in this book builds from sensori-motor operations into *conceptual embodiment* focusing on the properties of objects and *operational embodiment* focusing on operations, using language to describe and define more subtle forms of reasoning." (my italics)

The basic Van Hiele level uses the world as external, where the world supports human memory by providing the input to the senses, so that the mind can start collecting, memorizing, categorizing and so on. For didactics it is important to be aware what can be done at this level. In fact, in various other places Tall (2013) uses "embodied" in that external sense in various places, for example in tracing a curve by the fingertops to understand "tangent". If embodiment is external, then terms as "conceptual embodiment" and "embodied formal" are contradictions in terms.

If Tall intends embodiment to be internal (representation in neurons and chemicals) then he would lose contact with what Van Hiele proposes for the base level, and the didactics there. This kind of embodiment also becomes rather vague, without an operationalisation that a neopositivist would require - but which brain researchers might look for, such as Dehaene on numbers - and without much help for didactics that doesn't use brainscanners.

Tall is struggling with his terms overall. On p133 he refers to "thinkable concepts" but please explain "unthinkable concepts". On p425 he attributes his notion of "crystallization" / "crystalline" to a discussion with Anna Sfard, and on p 429 to a discussion with Koichu and Whiteley. It appears to be a rephrasing of Piaget's "encapsulation" and Van Hiele's compacting of various properties to a more unified concept. One better rejects that term "crystalline" since, taking a neopositivist stand again (though not in principle), it suggests more than there is. Are George's "crystalline" concepts the same as Harry's ? The phrase adds nothing, while of course the phrase suggests that Tall adds something to the discussion which on this point is not the case, as he only rediscovers what Piaget enlightened and what the Van Hieles systemized with proper empirics. (Van Hiele regretted the lack of empirical testing on his theory, but it remains empirical, and the evidence is supportive.) [See Colignatus (2015a) again on the "definition & reality methodology".]

Finally, Tall suggests the notions of *structural, operational and formal* abstraction, thus again relating to the two brain hemispheres and the prefrontal cortex, and relating to the new columns introduced with the new diagram / table for the Van Hiele level structure. Tall is careful enough to say that these types of abstraction may be difficult to distinguish. My problem with this is: (a) this may distract attention from the fact that the Van Hiele levels already concern abstraction, (b) this might come with the suggestion of something new but there isn't anything new, except for the words, (c) the distinction in types of abstraction creates an illusion of exactness.

Tall merely introduces new words to describe what happens in the Van Hiele levels. Van Hiele recognised that there can be a level shift in the step from arithmetic to proper algebra. He already knew that this was different from a level shift in geometry. What is the use to label the first as "operational" and the latter as "structural" abstraction ? And, confusingly, doesn't algebra contain some "structure" too ? Or doesn't geometry contain "operations" too (like drawing a circle with a compass) ? The "issue of distinguishing kinds of abstraction" is only created by the inclusion of the visual / aural columns, but this should not distract us. Overall, the use of new words may also be a matter of taste. The key point remains that there should not be the suggestion of something new on content.

Van Hiele rejected the vague word "abstraction". But once the Van Hiele levels in insight in mathematics are properly understood, I agree with Tall, and already proposed independently, that it can be advantageous to redefine the issue in terms of "abstraction" anyway, namely given that thought is essentially abstract and given that mathematics [deals with] abstraction.

## *Tall's view on Van Hiele and Freudenthal (2014)*

Setting the record straight w.r.t. Van Hiele as presented by Tall also causes a look at the influence of Freudenthal. If Tall did not fully appreciate the work by Van Hiele, might there have been an influence from what Tall read from Freudenthal's description of that same work ? Thus, while we focus in [this part of the 2014] article on VH & T, there is also F, generating the relationships VH & F and T & F. If the relationship VH & F was fair and F reported correctly, then this report would not have been a cause for T's long misunderstanding of VH. But if something went wrong in [VH & F] and the reporting was biased, then Tall might be wrongfooted, as the rest of the world.

In Freudenthal and Van Hiele we have two Dutch researchers with some impact on mathematics education on the world stage. It has some interest, both on content and history, what an outsider like professor Tall observes on this, and how this affected and affects his appreciation of their works. Note that the commission on math instruction of the international mathematics union (IMU-ICMI) has a "Hans Freudenthal Award" rather than a "Piaget - Van Hiele[s] Award". We may somewhat infer that Freudenthal's international influence seems greater. An international reappreciation of Van Hiele might cause an reconsideration though.

Let us thus focus now on Tall's perception and presentation of the Van Hiele and Freudenthal combination. As said, Tall is a foreigner and outsider to this, while the insider Dutch have the advantage of additional personal information, documents and e.g. newspaper articles in Dutch, but perhaps the disadvantage of missing the bigger picture. My own position comes with the advantage of distance in time, as I came to teaching math at highschool only in 2007. I have the (dis-) advantages of being Dutch and a foreign exchange highschool student year in California 1972-1973. Also my first education was a degree in econometrics in 1982 and my degree in teaching mathematics came later in 2008. My background in empirics differs from a first training in abstraction only, as happens with mathematicians.

As William Thurston (1990, 2005) and Hung-Hsi Wu (see Leong (2012)) complain for three decades about the dismal state of math education in the USA, one should hope that there are independent factors at work in the USA itself, but the influence from Freudenthal with his advocacy of "realistic mathematics education" (RME) should not be regarded as negligible. When Thurston submitted his (1990) text to arXiv, he added this comment (2005):

> "This essay, originally published in the Sept 1990 Notices of the AMS, discusses problems of our mathematical education system that often stem from widespread misconceptions by well-meaning people of the process of learning mathematics. The essay also discusses ideas for fixing some of the problems. Most of what I wrote in 1990 remains equally applicable today."

A key document to understand the Van Hiele - Freudenthal combination is the interview in Dutch with Van Hiele by Alberts & Kaenders (2005), in the mathematics journal NAW of the Dutch Royal Mathematical Society (KWG). There Pierre van Hiele says: *"Trouwens, Freudenthal heeft mij later nogal eens een hak gezet, jongens."* (p247). This is Dutch idiom. Google Translate July 26 2014 gives literally: "Besides, Freudenthal has often put a heel to me later, guys." My proposed free but clearer translation, taking account the rest of the interview, is: "Besides, Freudenthal has later frequently sabotaged my work, guys." It is important to add that Van Hiele remained polite, as a fine math teacher would do. [2015: See above for more translation.]

Broekman & Verhoef (2012:123), in their short biography after the decease of Van Hiele in 2010, confirm that Van Hiele would have appreciated a university research position but wasn't offered one, and thus remained a highschool teacher all his life. The reason is not that such positions could not have been made available, or even a professorship. Broekman & Verhoef describe differences of opinion between Freudenthal ("reality vs model") and Van Hiele ("concrete vs abstract") but they do not mention the crucial distinction that Van Hiele had an empirical attitude while Freudenthal remained locked in mathematical abstraction (with a virtual notion of "reality").

On my weblog, I have concluded that Freudenthal's "realistic mathematics education" (RME) (i) partly doesn't work and (ii) that the part that works was mostly taken from Van Hiele. Freudenthal's "guided reinvention" is Van Hiele's method of providing the students with the relevant materials and instructions so that they can advance in the levels of understanding. In World War 2 in 1940-1945 Freudenthal was in hiding, taught his children arithmetic, relied on real world examples, and wrote a notebook on this. Van Hiele provided theoretical justification for the concrete vs abstract distinction. The name RME essentially refers to the basic Van Hiele level. Though Freudenthal will have been inspired by his own experience, he must have realised that the scientific basis was provided by Van Hiele, and thus his choice of the name "realistic mathematics education" amounted to some appropriation and distortion of the Van Hiele result. While Freudenthal at first referred to Van Hiele he later tended to refer to "his own publications", which had the effect that Van Hiele wasn't openly referred to. Colignatus (2014a) concludes that RME is a fraud. [14] These conclusions on Freudenthal himself and RME directly affect our understanding of how Tall and the rest of the world were disinformed, which helps our understanding of Tall's perception of Van Hiele's work. My weblog short discussion and conclusion clearly only touch the surface of the problem. One can only hope that funds will become available both to analyse the errors of the decades since 1957, and to develop ways to repair those.

To put the issue in context: Colignatus (2014b) concludes that there is a serious issue with scientific integrity in the Dutch research in mathematics education, starting with Freudenthal but now with a dysfunctional "Freudenthal Institute" with a loyality complex, and with also failing supervision by the Dutch Royal Academy of Sciences (KNAW). The same disclaimer w.r.t. a weblog article applies. Those conclusions on the currently dysfunctional "Freudenthal Institute" are of less direct relevance for our topic of discussion [in this section] of getting the facts right on Van Hiele & Tall. However, it is useful to mention them since they put the issue in context. For example, professor Tall indicated (in an email conversation) that it was hard to find English sources of Van Hiele's work. Part of the explanation is not only in the observable sabotage by Freudenthal himself but also in the continued similar dysfunction by said institute, that should have been able to recognise the importance of Van Hiele's work and help make translations available.

---

[14] This particular weblog entry also makes fun of the method of "history writing" by Amir Alexander, and thus may require some decoding. The statement in July 6 2014 remains careful: "It might be that Freudenthal is no real crook but never waste a good story." The Addendum of July 15 however leaves little to guess: "Thus, Van Hiele was aware of the portent of his theory, contrary to what David Tall suggests. So much more of a pity that Freudenthal sabotaged and appropriated it."

After I wrote that analysis on my weblog (2014a), a Dutch reader drew my attention to the thesis in Dutch by La Bastide-Van Gemert (2006) (further LB-VG) that I had not seen yet at that time. This thesis is explictly on Hans Freudenthal on the didactics of mathematics, and her chapter 7 discusses Freudenthal and the Van Hiele levels. At first I was inclined to neglect this hint and study since it seemed that it did not pertain to this present paper on Van Hiele & Tall. But checking it, it appeared quite relevant. I see my weblog analysis confirmed. Above, below and in the **Appendix B** I provide some English translations of some passages. LB-VG describes on p195 that Freudenthal's work on math education was rather bland and traditional before 1957 and only gains content after the work by the Van Hieles. She decribes how Freudenthal first refers to the Van Hieles but then introduced new terms like "guided re-invention" and "anti-didactic inversion", and subsequently moved Van Hiele into the background and started advocating those own terms. However, LB-VG apparently falls for the suggestion that Freudenthal really contributed something new, and not just new phrases and misunderstandings. Given this (wrong) perspective, she does not discuss the scientific integrity problem of misrepresentation, appropriation of work and the withholding of proper reference. With the proper perspective her thesis however provides corroboration.

Above, I already quoted some evidence from her thesis concerning Pierre van Hiele's early awareness of the portent of his theory. There are some more particulars on Van Hiele and Freudenthal that might distract here and thus are put into **Appendix B**. A conclusion there is: While Freudenthal took key parts from Van Hiele's theory, he also inserted his own phraseology, with such consequence that Tall apparently had difficulty recognising Freudenthal's texts as Van Hiele's theory in (distorted) disguise, so that Tall could embark on his own path to re-invent Van Hiele's theory.

**David Tall and RME and Dutch language (2014)**

How does professor Tall deal with this situation of which he has been a foreign observer for all these decades ? Tall (2013:414-415) has a short text about 'realistic mathematics education' (RME). Tall (2013) refers to Van den Heuvel-Panhuizen (1998) of the "Freudenthal Institute", henceforth VdH-P. It is useful here to also refer to an early review by Tall (1977) of Freudenthal's book *"Mathematics as an educational task"*. (Note that I regard mathematics as abstract and education as empirical, so that Freudenthal's book title reads to me as a contradiction in terms.) Thus, Tall (2013:414-415) on RME:

> "The Dutch project for 'realistic mathematics education' was introduced to build on the learner's experience and to replace an earlier mechanistic system of teaching routine procedures. [footnote reference to VdH-P] It provides the child with a realistic context in which to make sense of ideas that are often performed in a practical situation. Yet, as time passed, it was found that, at university level in the Netherlands, remedial classes needed to be introduced because more students lacked the necessary skills for advanced work in mathematics and its applications." [footnote reference to "Information supplied by my colleague Nellie Verhoef, based on articles in Dutch: (....)]

> The Dutch translation of the verb 'to imagine' is 'zich *realise*ren', emphasizing that what matters is not the real-world context but the realization of *the reality in the student's mind,* which Wilensky expressed as the personal quality of the mental relationship with the object under consideration. [footnote reference]
>
> The three-world framework not only sees practical mathematics related to real-world problems, but it also offers a theoretical framework to realize ideas in a conceptual embodiment that transcend specific examples and blend with flexible operational symbolism."

We can deconstruct this quote on the points of (i) translation, (ii) the rewriting of history, and (iii) textual gibberish. Incidently, Van Hiele (1957) distinguishes between students who rely on algorithms and students who are able to recognize structure. Apparently the reference to "mechanistic" by VdH-P concerns an emphasis in teaching upon algorithms rather than insight. But the words algorithm, structure and insight are avoided.

Firstly, Tall is erroneous on the Dutch translation:

- The translation of English "to imagine" to Dutch is "zich verbeelden, zich inbeelden, zich voorstellen" (as alternatives with different shades of meanings) and *not* "zich realiseren". I checked that Google Translate July 26 2014 had it good.
- The translation of Dutch "zich realiseren" (note the "zich"= "oneself") is "to become aware", as "He realised that the train would depart without him" or "She realised that a triangle with two equal sides also has two equal angles". Or see above quote of Tall (2013:430) how he realized only at the death of Van Hiele how important his contribution had been (except that Tall regards this as his own discovery and not something that Van Hiele had already been aware of and explaining about).
- It is not *making-real* as in "He realised his plan to reach the top of the mountain". In Dutch "Hij realiseerde zijn plan" has such meaning. The difference comes from "zich" (oneself) between "realiseren" en "zich realiseren".

Tall's erroneous translation of "to imagine" to Dutch apparently is based upon the erroneous translation provided by VdH-P (1998). This VdH-P text is problematic in various respects. It is useful to realise that there already was quite some criticism on RME in 1998, so that VdH-P's text also has a quality of defense against that criticism. It also rewrites history and presents a curious view on education. Apparently unwittingly, Tall is dragged along in this. VdH-P (1998):

> "The present form of RME is mostly determined by Freudenthal's (1977) view about mathematics. According to him, mathematics must be connected to reality, (...) It must be admitted, the name "Realistic Mathematics Education" is somewhat confusing in this respect. The reason, however, why the Dutch reform of mathematics education was called "realistic" is not just the connection with the real-world, but is related to the emphasis that RME puts on offering the students problem situations which they can imagine. The Dutch translation of the verb "to imagine" is "zich REALISEren." It is this emphasis on making something real in your mind, that gave RME its name. For the problems to be presented to the students this means that the context can be a real-world context but this is not always necessary. The fantasy world of fairy tales and even the formal world of

> mathematics can be very suitable contexts for a problem, as long as they are real in the student's mind."

Deconstructing this:

(1) First of all, the true origin of the "realistic" in RME derives from Freudenthal's real-world linkage and not from the "zich realiseren" (become aware). LB-VG (2006:194) confirms that Freudenthal in persecuted hiding during the war years 1940-1945 already wrote a didactic text for teaching arithmetic to his children, in which real life issues are prominent. This also relates to the switch from abstract Sputnik "New Math" to applied mathematics around 1970. But the theoretical justification for the concrete vs abstract opposition and the start from mere intuition only derives from the Van Hiele theory. Van Hiele better explains than Freudenthal that teaching should start from what students understand, and that teaching is guiding them towards what they don't understand yet. The overall conclusion is that Freudenthal took the RME name from the Van Hiele basic level, as that had received theoretical justification.

(2) Secondly, the "mostly" and "not just" are a rewriting of this history into another interpretation, in which the "zich realiseren" is plugged in, taking advantage of the flexibility of language. There seem to be earlier occurrences of "zich realiseren" (possibly Gravemeijer 1994) before this present use by VdH-P, but this should not distract from the effort at rewriting history w.r.t. point (1) after 1957.

(3) Thirdly, while the Dutch term "zich realiseren" has an etymological root to "reality", the subsequent explanation by VdH-P of "zich realiseren" should be about "to grow aware" (which is the proper translation). However, we see that the explanation that she gives is about the ability to understand (imagine) the situation under discussion. While history is rewritten, this new interpretation of "realistic" has a different meaning in Dutch (to grow aware) than in English (to be able to imagine). A wrong translation is deliberately used to suggest that the meaning would be the same in both languages.

(4) On content, I fail to understand why it would make a difference whether a student sees a fantasy **as** *the fantasy that it is* **or** *experiences it 'as real'*, if you want to link to their existing stock of experience and mental abstraction to start doing mathematics. The suggested didactic condition "as long as they are real in the student's mind" is unwarranted. (For example, if one presents a fantasy cartoon image of a blue car and a photograph of a red car, and starts a discussion whether one can do 1 car + 1 car = 2 cars, we may suppose that part of the discussion would be about differences between cartoons and photographs, and that you cannot add images to real cars on the parking lot. But you can count cars on the parking lot and add their numbers and use images to represent them, even fantasy cartoons and "realistic" photographs.)

Subsequently, Tall's phrase "realize ideas in a conceptual embodiment" is gibberish. The "realization of ideas" is ambiguous, as said. Is it growing aware or is it making-real ? Students can draw lines and circles, to approximate abstract ideas, but we cannot assume that they can turn cartoons of cars into actual cars on the parking lot. To grow aware of ideas is excellent, but this requires the Van Hiele teaching method and not the distorted versions. The "conceptual embodiment" is, as said, either a contradiction in terms, or it refers to a constellation of neurons and biochemicals, but then yet lacks an operational meaning useful for teaching.

Taking Tall's quote as a whole, he seems to suggest that his confused approach would be helpful for the Dutch situation that pre-university teaching has gone haywire and requires remedial teaching at university. His frame of reference prevents him from observing that the error lies in RME and in the work of the very VdH-P whom he quotes. It prevents him from "realizing" that he himself is also off-track in his book and in joining up with a failing Dutch constellation.

What really happens here is that professor Tall apparently wants to link up to the Dutch situation and RME theories, perhaps wanting that his theory is also accepted at the dysfunctional "Freudenthal Institute", without having adequate understanding about the local maltreatment of Van Hiele's work. The information that he relied on by VdH-P was wrong. The information that he quoted as receiving from Verhoef seems incomplete, and last week [August 2014] I got evidence that criticism is possible on subsequent issues. Who in Holland will defend Pierre van Hiele, with the indicated loyality complex at the "Freudenthal Institute" ? I informed Tall about my books Colignatus (2009) [now with second edition 2015] and (2011a). Professor Tall and I actually met at a conference in Holland on June 24 2010 and spoke a bit longer. In an email of 21 May 2012 he mentioned that he was in the autumn of his life and wished to get his book finished, and didn't have time to look at my books. It is a regrettable paradox that the amount of time taken to write is very much more than the amount of time needed to see that one should write something else and much shorter. Still, one can only respect a person in the autumn of his life. In a way I consider it very useful that Tall has taken stock of his work, since it shows both a misunderstanding of Van Hiele over most of Tall's life and a recognition of the importance at that autumn. Hopefully Tall continues to think about math education, hopefully also in an essential rewrite.

### *Conclusion on Van Hiele and Tall (2014)*

It is somewhat enlightening to conclude that professor Tall as researcher and teacher in the education of mathematics seems to have little experience in the use of Van Hiele methods in actual educational practice. In itself it is relevant that Tall recognises the importance of Van Hiele's work at this late stage, in retirement. It remains also a phenomenon to be explained that such an important theory recieves such recognition by Tall only at such a late stage. For this, the details of the Van Hiele - Freudenthal combination are relevant, of which researchers in education in mathematics do not appear to be aware about in general.

David Tall grew aware to a much larger extent in 2010 of the Van Hiele niveaux of understanding of mathematics. He also thought that Van Hiele (1957, 1959, 1973, 1986, 2002) saw only limited application. Tall now claims that it was a creative insight on his part to extend those levels to wider applicability. Perhaps it was, given his misunderstanding of Van Hiele. Perhaps it was only a recollection of something read or heard but forgotten and surfacing in different form. Whatever this be, the claim however does no justice to Pierre van Hiele who already asserted that wider applicability, also for other disciplines than mathematics, in 1957. Tall's claim may block researchers in education in general and the education of mathematics in particular from considering the wealth in Van Hiele's work. We owe Pierre and Dieke van Hiele and our students to get the facts right.

## *More information in 2015 on Van Hiele and Freudenthal*

This article in 2014 focussed on Tall given his 2013 book. The findings on Freudenthal were supplementary. Freudenthal apparently got greater weight in the subsequent discussion because his misconceptions on "realistic mathematics education" (RME) have been around for longer with also more impact on society. This version of the article in 2015 inserts Freudenthal in the title. For the line of reasoning in 2015 the reader should now read **Appendix B**. It is still kept in its original form and position, and not included in the body of the text, so that the reader can reconstruct the situation in 2014 and subsequent discussion. Also, **Appendix B** contains my 2014 translations of the 2006 LB-VG thesis, while there now is an English translation approved by LB-VG herself.

There now is also a review by Selden (2015) in which an expert with no ties to Holland considers the LB-VG thesis in its English translation. Selden (2015) however overlooks the inconsistency in the thesis.

## *Van Hieles 1957, full paragraph versus partial quote*

The Euclides digital source also shows that LB-VG gives only a partial quote. (With the reference we could have gone to the library in 2014.) On page 202 of her thesis (Dutch 2006) LB-VG states, my translation and emphasis:

> "Now the Van Hieles had thought about it themselves as well to apply the theory of levels to other subjects **in mathematics education.** Already in 1957 the Van Hieles indicated, in an article in Euclides about the phenomenology of education that gives an introduction to geometry, **not to exclude that possibility**: [quote and footnote 82]" [15]

As already concluded in **Appendix B** in 2014, this judgement is inconsistent with what LB-VG quoted, on the very same page 202, from the 1957 article: "**The approach presented here for geometry namely can be used also for other fields of knowledge (disciplines).**"

In 2014 I had only LB-VG's partial quote of the 1957 Euclides article. It has been fully quoted above. LB-VG leaves out these lines – and when you drop these then it is easier to suggest that the Van Hieles weren't quite aware of what they were doing:

> "For students, who have participated once in this approach, it will be easier to recognize the limitations than for those students, who have been forced to accept the logical-deductive system as a ready-made given. Thus we are dealing here with a formative value (Bildung), that can be acquired by the education in the introduction into geometry."

---

[15] Dutch: "Nu hadden de Van Hieles er zelf ook wel aan gedacht de niveautheorie ook op andere onderwerpen uit het wiskundeonderwijs toe te passen. Al in 1957 gaven de Van Hieles in een artikel over de fenomenologie van het aanvankelijk meetkundeonderwijs in Euclidesaan die mogelijkheid niet uit te sluiten:"

With the full quote available, LB-VG's judgement on the 1957 Euclides article becomes even more curious. As I write in Colignatus (2015c):

- The Van Hieles do **not** limit this general applicability to mathematics education only. They speak about *other fields of knowledge*. It is LB-VG who puts it into the box of mathematics education only.
- It isn't "not excluded" but **emphasized.**
- The general claim is in the theses (ceremony July 5 1957) under supervision of Freudenthal and not just the article (October 1 1957).
- Why not quote the full paragraph ? It would show why the Van Hieles select geometry also for its ability by excellence to teach this **general** lesson. Pierre van Hiele's thesis has the word **"demonstration"** in the title, to that the discussion of geometry is only intended to demonstrate the general applicability (in the same manner as demonstration is used in geometry itself).

Did LB-VG first cut out a quote and only then think about it ? Or, did she already have the frame of mind that it was Freudenthal who really found the general applicability, such that the selection of Van Hiele sources had to fit into that frame of mind ? Still, the part that she quoted does cause inconsistency in her thesis.

### *English translation in 2015 of the 2006 thesis by LB-VG*

This article is about getting the facts right: by stating those facts. There now arises an element of *trying to get* some facts. This paper is about Van Hiele, Tall and Freudenthal, and not about LB-VG and her interpretations, handling of sources and priorities w.r.t. the discovery of the inconsistency in the thesis. Some comments however are required, e.g. for readers who compare **Appendix B** (that has my translations of quotes from the Dutch 2006 thesis) with the new 2015 English translation of the full thesis.

On August 19 2014 I alerted La Bastide – Van Gemert (LB-VG) to the inconsistency in her Chapter 7, and referred to the 2014 version of this article on my website with the full discussion (also available on arxiv.org). I suggested that she had a rosy view on Freudenthal and likely wasn't aware of what was happening here. I expressed my hope that she would look into this. LB-VG replied the same day that she didn't have time for this. The two emails are in **Appendix C** (in Dutch). My translation of her answer:

> "Thanks for your email. Unfortunately I don't have the possibility now to respond on content. I wish you much luck with your article."

I presumed that her time was required for research in epidemiology at the University Medical Center Groningen (UMCG). Given her rejection I asked in 2014 professors Klaas van Berkel, Jan van Maanen and Martin Goedhart, all involved with the 2006 thesis, and still in academic positions, whether they could look into the inconsistency in the thesis. They all declined. Their neglect to look into it is problematic. The inconsistency is not something hidden, or something complex that requires long thought. Professors Van Berkel, Van Maanen and Goedhart should have spotted it immediately in 2006 – but we don't know the phases of the thesis. They should have recognised it immediately once I pointed it out in 2014 and asked them to look into it.

I was very surprised to discover the 2015 English translation of her thesis, see Colignatus (2015h).

> **Problem:** LB-VG didn't have time to look into this inconsistency but did have time to have it translated ?

A lawyer might argue that the thesis was established in 2006 and that this concerns a mere translation into English of that same work. A scientist, however, who has been alerted to an inconsistency doesn't merely translate. One can translate a text *as it is*: but one will also include a note about the discovery of the inconsistency. It is even possible that the inconsistency causes that a translation becomes rather useless. While the thesis portrays Freudenthal as some hero – who supposedly discovered what Van Hiele overlooked – we however found intellectual theft and fraud. Thus it seems that LB-VG knowingly dispatches wrong information into the English speaking world. Since there may be circumstances that I am not aware of – perhaps she is terminally ill and this English translation was her final wish ? – I have submitted the problem to the integrity officers of UMCG. Until I have information to the contrary I tend to feel misinformed about that lack of time.

Colignatus (2015c): It wasn't just my question in 2014 that might have alerted her. The Alberts & Kaenders interview of 2005 with Pierre van Hiele's protest was in the year preceding the 2006 thesis. It is not in her list of references. We might accept that the interview was published too shortly before the thesis to affect it, but the period 2006-2014 would normally have allowed a review of the argument in the light of Van Hiele's protest. His decease in 2010 didn't go unnoticed. It is interesting in itself to see that also a reviewer like Danny Beckers (2007) apparently neither read that interview, for he characterises Van Hiele as a "friend" of Freudenthal. There are strong "frames" that appear to cause people to overlook Freudenthal's fraud. (One should hope that Beckers now also translates his review into English, and takes along the new information.)

**Key point 1: Independence**

What is important for the reader to know: *I have not looked at the 2015 LB-VG English translation*. I presume that it contains the same arguments as the Dutch original, and as translated in **Appendix B**. (If some phrasings have changed, that would cause also a new Dutch edition, then I hope that someone informs me.)

I see this expectation confirmed by the review by professor Annie Selden (2015) on the website of the *Mathematical Association of America* (MAA). Selden praises the book, and *doesn't see the inconsistency*, see Colignatus (2015i). With my emphasis:

> Selden states: "Chapter 6 covers the period from 1950 to 1957 when Freudenthal's national and international reputation as a mathematics educator grew enormously. Also, towards the end of the period, his mathematical-didactical ideas were greatly influenced by the pedagogical dissertation studies of Pierre and Dina van Hiele **on geometry.** Their work, and its influence on Freudenthal who was their dissertation advisor, is further discussed and analyzed in Chapter 7."

Colignatus (2015b) discusses Freudenthal (1948) on didactics. Colignatus (2015h): In the period before 1957 Freudenthal's ideas on education are rather bland. LB-VG describes how they grow into RME only after 1957, when he has the theses by the Van Hieles. Again, there is the suggestion that Pierre van Hiele looked only at geometry, while he stated the general relevance, and used geometry only for demonstration (with a wink reference to the role of demonstration in geometry).

> Selden states: "Of special interest to mathematics education researchers who use realistic mathematics education (RME) as their theoretical framework is Section 7.4 titled, *"Freudenthal and the theory of the van Hieles: From 'level theory' to 'guided re-invention'"*. According to the author, it was during this time period that Freudenthal introduced the ideas of "guided re-invention" and the "anti-didactical inversion". These terms "did not come out of the blue. … [B]oth concepts were already mentioned before in more guarded terms. But it is the first time that Freudenthal mentioned and defined them explicitly." (p. 195)."

Colignatus (2015h): This is however where intellectual theft takes place. It is amazing that Selden doesn't observe it, but, she might not know the work by the Van Hiele's so well. The key question is whether the notions of "guided re-invention" and "anti-didactical inversion" are deep and special. If they would be, then Freudenthal could claim major discoveries. In fact, they turn out to be simplistic rephrasings of what Van Hiele already described. Van Hiele was interested in insight, and transitions to higher levels of insight. Now, isn't invention the phenomenon of arriving at more insight ? It is basically just another word. The same holds for the Van Hiele process from concrete to abstract, that is opposite to Euclid's *Elements* that starts with proofs. It is a bland rephrasing. Education in 1957 didn't have a refined taxonomy such that the lawyers of the Lesson Study inquisition could haggle about student *A* having a Van Hiele level transition and student *B* having a Freudenthal guided re-invention, with numbers to show that Freudenthal made the more relevant discovery. The conclusion is that these are just rephrasings, and that Freudenthal could, once he had his own terms and publications, refer to his own work rather than Van Hiele.

Colignatus (2015h): PM. If you like to think about the difference between cars and ideas, then there is this argument. You might suggest that new ideas are always your own. Thus Freudenthal's new phrases would still be something of his own, and he could always claim credit for them. For cars, this would mean that if the robber puts a new paint on your car, he can keep it. It is an interesting suggestion. It would also hold when the robber puts so much paint on the radiator and exhaust that the car would hardly run, like RME hardly works. Thus, think about it. Your dear car, stolen and turned into a wreck, with the robber dancing and prancing atop.

Incidently, LB-VG (206:201) pointedly refers to Freudenthal's connection to Brouwer's Intuitionism here, in which mathematics is constructed. See Ernst Snapper (1979) for an excellent exposition on the foundations of mathematics. See Colignatus (2015l) on foundations and degrees in constructivism.

**Key point 2: Retractions**

The 2015 English translation of the 2006 LB-VG thesis must be retracted, and at least must be extended with a warning on the inconsistency. [16] The Selden review should be retracted too as it is based upon a disinformative book, whence there is involvement by the MAA Reviews editor and the MAA executive, see Colignatus (2015i).

**Key point 3: Groningen**

This issue of integrity at UMCG and RUG appears to become complex. (a) UMCG and RUG have complex governance on scientific integrity, making it difficult for outsiders to determine who is responsible for what. (b) Two UMCG integrity officers decided not to ask LB-VG for her copies of the email exchange in 2014-2015, and they "closed the case" with the argument of lack of information: which is a false argument when you didn't look for the information yourself. I have asked UMCG / RUG to replace these officers. I suppose that nothing will happen in the mean time, which isn't kind to LB-VG and Selden and libraries and people who buy the book at Springer. (c) One should also look at the conduct of the three professors involved with the 2006 thesis who declined in 2014 to look into the inconsistency when it was pointed out. (d) Presently, I regard the UMCG / RUG integrity situation as a disaster zone. Emails are available here. [17]

## *Prospect on neoclassical mathematics*

This article is not about some old men. A key notion is the distinction between mathematics as dealing with *abstraction* and education as an *empirical* issue. Freudenthal had been trained for pure mathematics and presumably it fitted him. Van Hiele got his degree in mathematics but his heart was in education. Van Hiele relied upon his observations in practice. He was for the use of statistical methods but mentioned the obstacles that we are all familiar with. Freudenthal threw mathematics education research back to pre-scientific times by relying on his a priori perception of what he thought what would work. Apparently he used the work by Van Hiele to craft stories that others might buy. He got some momentum because of the New Math disaster by other abstract thinking mathematicians. The Freudenthal disaster ought to be a lesson for mathematicians: to stop meddling in mathematics education.

Progress is possible when such lessons are learned. It will be useful to give an example, taken from Colignatus (2015k). Namely, consider mixed fractions. The classically correct expression 2 + ½ (two and a half) has a traditional notation 2½ that however is confusing for pupils, for it has the (two times a half) structure like in 2√2. Pupils don't control spacing in handwriting as tight as a typewriter does. When they read or write 2½ as 2 ½ then they might defend outcome 1.

---

[16] It is a bit curious to put the following consideration here: If there is a decent explanation for the priority given by LB-VG to the English translation instead of looking into the inconsistency, then integrity of science has not been breached, and then it would generally be best that she first looks at the problem and the implications, and then decides for herself how she re-evaluates the book and what would still be relevant in it: and writes a position paper on this. A three months furlough should be enough to write such a position paper. The University of Groningen (RUG) allowed the inconsistency in the thesis, and this comes with some responsibility. If LB-VG is in the position to re-evaluate the text and its consequences, RUG should provide for those three months such that she can.

[17] http://thomascool.eu/Papers/BHRM/Index.html

Consider the simplification of 2½ / 3⅓. Van der Plas (2008, 2009) reports that this kind of question hardly occurs in Dutch "realistic mathematics" textbooks. Traditional mathematicians like Hung-Hsi Wu in the USA want to see a lot of practice on this again. In handwriting by students: 2 ½ / 3 ⅓ = 1 / 3 ⅓ = 1 / 9. Van de Craats in Holland is aware of the issue and suggests using 5 / 2 but then loses the location on the number line.

Van Hiele (1973) (likely also in (1986)) proposed the abolition of fractions by using the notation of the inverse, i.e. with exponent -1. My suggestion is to use $H$ = -1 (pronounce: eta). This gives neoclassical: $2 + 2^H$. In this way the use of -1 can be avoided. Pupils who are still learning arithmetic and who see -1 might think that they must subtract something. Instead they can learn the rule that $x \, x^H = 1$ provided that $x \neq 0$. Later when powers and roots are introduced then they can see the numerical value for $H$. See **Appendix E** for an example and Colignatus (2014e) for a longer discussion.

Thus, we can create Table 2 with the different approaches. The classical and neoclassical notations are proper mathematics, while the traditional notation so much favoured by mathematicians like professor Wu is crooked, cumbersome, errorprone, counterproductive: so-called "mathematics". The use of $H$ is still only a suggestion, and empirical tests must show whether it really works better. Subconclusions are:

- Attention for empirics may result into didactic improvement.
- The distinction between general and particular is important. This finding on mixed fractions is a particular instance and no statement on professor Wu's work in general.
- However, see Colignatus (2009, 2015) for more cases for didactic improvement.
- It is also useful to observe that Holland is trying to recover from "realistic mathematics education" (RME), but that Freudenthal exported it to the USA, and that it now boomerangs back via the "21st century skills" supported by the OECD.

**Table 2: Mathematics versus "mathematics"**

|  | *Mathematics (empirics, engineers)* | *"Mathematics" (mathematicians)* |
|---|---|---|
| *New* | Neoclassical: $2 + 2^H$ | "21st century skills" (mostly old wine) |
| *Old* | Classical: 2 + ½ | Traditional: 2½ |
|  |  | "Realistic mathematics": try to hide this |

## *Conclusions*

Why would a professor of mathematics want to create his own theory of education, and tell a teacher of mathematics that his approach developed in practice doesn't work ?

When Freudenthal found his mathematical powers lessening, he offered himself the choice of doing history of mathematics or education of mathematics. His wife Suus Freudenthal-Lutter was involved in education, and there were other influences as well. Remarkable is that the university in that period allowed such an easy switch, while it would have been more logical to ask Pierre van Hiele to take a professorship in education of mathematics. The history of these events isn't so complex.

The conclusions are in the abstract and need not be repeated here.

*Thomas Colignatus is the science name of Thomas Cool, econometrician (Groningen 1982) and teacher of mathematics (Leiden 2008), Scheveningen, Holland.*

*I thank professor David Tall for graciously even though critically informing me and providing details so that I could better understand his analysis and position. I have offered to write a joint paper on this issue so that others would not be in doubt concerning his reaction to this information that apparently is new to him. Given his initial rejection of this suggestion I deem it better to clearly state that information.*

*I also thank professor Jan Bergstra (Amsterdam, KNAW) for providing critical comments and for drawing my attention to La Bastide - Van Gemert (2006). As the latter thesis focuses on Freudenthal I at first considered it not relevant for the focus in 2014 on Van Hiele and Tall, but it is another indication of the interrelatedness of things that it appeared to contain very relevant evidence. In 2015 I thank Bergstra for asking why the critique w.r.t. Freudenthal cannot be restricted to saying that he should have referred more. Bergstra may evaluate the evidence differently – but he has been trained on abstract mathematics.*

*Appendix A: Basic data and problems of translation (2014)*

(1) Pierre van Hiele's thesis advisor was H. Freudenthal while Dieke van Hiele's thesis advisor was M.J. Langeveld. The main source is http://dap.library.uu.nl. These data are also in http://genealogy.math.ndsu.nodak.edu/id.php?id=102372 and id=144944. Some other sources give conflicting information:
1. The "Freudenthal Institute" has Langeveld as supervisor for Pierre but this is only correct for his role as second supervisor:
http://www.fisme.science.uu.nl/wiki/index.php/Pierre_van_Hiele
2. Broekman & Verhoef (2012:123) state: "Van Hiele himself mentioned at more occasions, in private and not in public, that he interpreted that Freudenthal's choice to be the first supervisor of Dieke, which implied that he himself 'got' Langeveld as first supervisor, also indicated that Freudenthal rejected his more theoretical work (on the psychology of cognition and learning)." [18] This remark by Pierre thus would refer to the early period before the final thesis advisors were allocated.

(2) Professor Tall alerted me to Fuys et al. (1984) *"English Translation of Selected Writings of Dina van Hiele-Geldof and Pierre M. van Hiele"* downloadable at ERIC. Page 8 states that the translations have been accepted by Pierre van Hiele.

This source puts an emphasis on geometry, possibly stimulating the confusion for some authors that the Van Hieles might think that the levels applied to geometry only. Indeed, the translations were produced under the project title: *"An Investigation of the van Hiele Model of Thinking in Geometry among Adolescents".* Such a positioning might run the risk of excluding the wider scope. However, Dieke's chapter XIV *"Further analysis and foundation of the didactics"* looks at general didactics, in particular p202 when she compares with her own learning process in didactics. Importantly, the "tenets" (p220-2) start with a general scope (e.g. "I. In order to be able to arrive at an efficient study of a certain subject (...)".

Secondly, accepting a translation is another issue than finding a translation that reduces confusion. A key point is the translation of the title of Pierre van Hiele's thesis:
- My translation: *"The Issue of Insight, Demonstrated with the Insight of School Children in the Subject Matter of Geometry."* [19]
- The ERIC p8 translation, apparently accepted by Van Hiele: *"The Problems of Insight in Connection with School Children's Insight into the Subject Matter of Geometry".*

I have these considerations for my translation: (a) Keeping the original "demonstrated" in the title warrants that geometry does not only provide an existence proof for the levels but also forms only an example. (b) The reference to the role of demonstration in geometry itself must have been quite deliberate. The thesis discusses that mathematics is about *proof,* after all, and not just the execution of algorithms. (c) Van Hiele (1957) and (1959) emphasizes the triad teacher-student-subject. Thus, while the theory of levels is general, the subject matter exerts a relevant influence for particular educational situations.

---

[18] My translation of: "Zelf heeft Van Hiele in besloten kring meerdere malen het idee geopperd dat de keuze van Freudenthal om eerste promotor te zijn van Dieke, waardoor hijzelf Langeveld als eerste promotor 'kreeg', door hem altijd beschouwd is als een afwijzing door Freudenthal van zijn meer theoretisch (denk/leerpsychologisch) werk."
[19] My translation of: *"De Problematiek van het Inzicht, Gedemonstreerd aan het Inzicht van Schoolkinderen in Meetkunde-leerstof".* (The ERIC text on p258 suffers from typing errors.)

PM. Some smaller comments on Van Hiele (1959): (a) Van Hiele speaks about five levels and at first only presents four, but eventually the fifth appears on page 254 (using the ERIC page count). (b) It remains awkward that Van Hiele gave the label 0 to the lowest level. The ranking words "first" to "fifth" tend to become ambiguous, as they apparently also might function as adjectives (associating "the third level" with "level 3"). (c) The article is a translation from French, in which the word "intuition" might have been used for "insight" ? Overall, it might be better to call the base level the "intuitive level", where insights are still unguided (while there may of course be trained intuition at higher levels).

## *Appendix B: Additional information from the LB-VG thesis (2014)*

La Bastide-Van Gemert (2006) (further LB-VG) gives additional evidence on Freudenthal and the Van Hiele level theory. Her text apparently requires some deconstruction however since she appears to have a rather rosy view on Freudenthal's performance.

The body of the text above contains her quote of the *Euclides* 1957 article with the statement by the Van Hieles that the theory of levels applies to other disciplines too.

LB-VG then arrives at this curious statement:

> p197-198: "In that manner Freudenthal described the theory of levels by which the direct link with geometry, essential and explicitly relevant in the work by the Van Hieles, had disappeared. It seemed that Freudenthal by this abstraction started to see the theory of levels as independent of the context (geometry education) under which he had learned about it: a new level had been reached ...." [20]

Comment: (a) The suggestion of a direct essential link with geometry is inconsistent with the earlier observation that geometry was only an example. (b) LB-VG suggests that the thesis supervisor Freudenthal would not have seen that. (c) There is the false suggestion that only Freudenthal made that step into abstraction toward general application.

> Just as curious, p198: "In that manner, step by step, Freudenthal gave his own interpretation of the theory of levels. Independent from the education of geometry from which the theory originated, he abstracted it into a method of logical analysis in the clarification of the (levels in the) educational topic of interest." [21]

Comment: Apparently she is not aware of the inconsistency, and would not quite understand what Van Hiele had achieved.

p199. Another point, of which Van Hiele will have been aware, but which apparently was also claimed by Freudenthal, and again by Tall (2013), chapter 14, but apparently new to LB-VG (2006):

> "In modern mathematics the mathematical systems, that have arisen by the organisation and ordering of the topic of interest (the relevant mathematical issues), became the subject of organisation themselves, via axiomatisation. We find a remarkable parallel between this remark [by Freudenthal] and

---

[20] My translation of: "Zo beschreef Freudenthal de niveautheorie op een manier waarbij de directe link met de meetkunde, essentieel en nadrukkelijk aanwezig in het werk van de Van Hieles, verdwenen was. Het leek erop dat Freudenthal door deze abstrahering de niveautheorie los begon te zien van de context (het meetkundeonderwijs) waarin hij ze leerde kennen: een nieuw niveau was bereikt...."

[21] My translation of: "Zo gaf Freudenthal de niveautheorie van de Van Hieles stap voor stap een eigen invulling. Los van het meetkundeonderwijs waar de theorie uit voort kwam, abstraheerde hij deze tot een werkwijze van logische analyse bij het inzichtelijk maken van de (niveaus van de) leerstof."

> Freudenthal's interpretation of the theory of levels: in mathematics itself there were, in this manner, transitions (definable by logic) to a higher level, comparable to the transitions between the levels such as there appeared to exist within the process of education." [22]

While the above indicates that LB-VG would not have been at home in the Van Hiele theory of levels, the following in her chapter 7 gives evidence on other forms of scientific malconduct by Freudenthal:

- p191: A 1957 newspaper article quotes Freudenthal giving a wrong description of the levels (namely: in mastery of routine) - but such may happen with newspapers.
- p204: "The theory of levels as such disappeared in Freudenthal's publications into the background and he based himself primarily upon his own ideas such as 'anti-didactic inversion' and 'guided re-invention' that for him were related (implicitly or not) to this theory of levels. He still mentioned the Van Hieles and their work in his articles but now only in passing. For Freudenthal the work by the Van Hieles had been promoted to basic knowledge." [23] Comment: In my analysis, LB-VG takes a rosy view on this. Freudenthal must have known that the Van Hiele levels were not well-known, and certainly not their claim on wider application than geometry only. Instead, Freudenthal reduced Van Hiele to geometry only and advanced his own phraseology as the proper approach in general.
- p182 gives a quote by Freudenthal in his autobiographic book p354, which is rather convoluted and lacks the clarity that one would expect from a mathematician: "The process of mathematisation that the Van Hieles were mostly involved with, was that of geometry, more exactly put: they were the first who interpreted the geometric learning process as a process of mathematisation (even though they did not use that term, and neither the term re-invention). In this manner Pierre discovered in the educational process, as Dieke described it, the levels of which I spoke earlier. I picked up that discovery - not unlikely the most important element in my own learning process of mathematics education." [24] Comment: Freudenthal thus suggests: (a) Pierre's insight is just seeing what Dieke described, so that she would be the real discoverer. (b) Freudenthal's words "mathematisation" and "re-invention" would be crucial to describe what happens in math education, otherwise you do not understand what math education is about, and it is only Freudenthal who provided this insight. (c) The Van Hieles wrote about geometry but were limited to this, so that it was Freudenthal himself who picked it up and provided the wider portent by means of his new words.
- p194, taking a quote from Freudenthal's autobiographic book p352: "What I learned from the Van Hieles I have reworked in my own manner - that is how things happen." [25] Comment: This is the veiled confession of appropriation. Freudenthal claims to be powerless and innocent of deliberate appropriation since "that is how things happen". Who however considers what that "reworking" involves sees only phraseology and lack of proper reference.
- p205: She discusses Freudenthal's use of "reflection" for the level transition, that Van Hiele (2002) protests about. As far as I understand this discussion, Freudenthal essentially merely provides introspectively, and without empirical support, the word "reflection" in relation to a level transition, as if only his new word is the valid approach, so that only he can be the inventor of transition via that proper word. However, the proper scientific approach would have been to describe what the Van Hiele theory and approach was, then define what the new idea of reflection would be, and provide the empiral evidence on that new insight (as the Van Hieles had provided empirical data for their method to achieve level transitions).

---

[22] My translation of: "In de moderne wiskunde werden de wiskundige systemen die ontstaan zijn door het organiseren en ordenen van het onderwerp (de betreffende wiskundestof) zélfonderwerp van organisatie, van axiomatisatie. Tussen deze opmerking en Freudenthals interpretatie van de niveautheorie is een frappante parallel te trekken: in de wiskunde was er op die manier sprake van (door logica definieerbare) sprongen naar een hoger niveau, vergelijkbaar met de sprongen tussen de niveaus zoals die er binnen het onderwijsproces bleken te zijn."

[23] My translation of: "De niveautheorie als zodanig verdween in Freudenthals publicaties naar de achtergrond en hij beriep zich voornamelijk op de voor hem (al dan niet impliciet) met deze theorie samenhangende ideeën als 'anti-didactische inversie' en 'geleide heruitvinding'. Hij noemde de Van Hieles en hun werk nog steeds in zijn artikelen, maar nu slechts en passant. Het werk van de Van Hieles was voor Freudenthal gepromoveerd tot basiskennis."

[24] My translation of: "Het mathematiseringsproces waar de Van Hieles zich vooral mee bezighielden, was dat van de meetkunde, preciezer gezegd: ze waren de eersten die het meetkundig leerproces als proces van mathematiseren interpreteerden (al gebruikten ze de term niet, evenmin als de term heruitvinding). Zodoende ontdekte Pierre in het onderwijs, zoals Dieke het beschreef, de niveaus waarvan ik eerder sprak. Ik pakte die ontdekking op – wellicht het belangrijkste element in mijn eigen wiskunde-onderwijskundig leerproces."

I might mention that LB-VG doesn't seem to be aware of Van Hiele's insight in the role of language. We could consider more points but it seems that the above suffices.

Another conclusion is: While Freudenthal took key parts from Van Hiele's theory, he also inserted his own phraseology, with such consequence that Tall apparently had difficulty recognising Freudenthal's texts as Van Hiele's theory in (distorted) disguise, so that Tall could embark on his own path to re-invent Van Hiele's theory.

---

[25] My translation of: "Wat ik van de Van Hieles leerde heb ik op mijn eigen wijze verwerkt – zo gaat dat nu eenmaal."

## *Appendix C: A single email in 2014 of TC to LB-VG, and a single reply*

My email of 2014 that asked about the inconsistency (Dutch original):

> Van: Thomas Cool / Thomas Colignatus
> Verzonden: 19 August 2014
> Aan: Bastide-van Gemert, S la
> Onderwerp: N.a.v. uw proefschrift, Hoofdstuk 7, Van Hiele niveaux
>
> Geachte dr. La Bastide – Van Gemert,
>
> Heeft u nog interesse in uw proefschrift, of bent u doorgegaan naar nieuwe terreinen, zoals de epidemiologie ?
>
> N.a.v. de situatie in het onderwijs in wiskunde en rekenen kwam ik ertoe ook te kijken naar de invloed van Hans Freudenthal.
>
> Relevant leken daartoe ook de herinneringen van David Tall, een Engelsman. Toen ik hem e.e.a. navroeg begon hij ook over Pierre van Hiele en zijn eigen jongste boek (2013).
>
> E.e.a. leidde tot dit artikel, beoogd voor een tijdschrift:
>
> http://thomascool.eu/Papers/Math/2014-07-27-VanHieleTallGettingTheFactsRight.pdf
>
> Uw Hoofdstuk 7 bespreek ik op p7-8. Mijn conclusie is dat u eigenlijk inconsistent bent, wanneer de Van Hieles in 1957 in Euclides al een algemene geldigheid voor hun theorie claimen, door u geciteerd, en u tegelijkertijd stelt dat Freudenthal dat pas aanbracht. M.i. heeft u dan een roze bril t.a.v. Freudenthal gehad, en niet doorgehad wat hier allemaal gebeurde. Per saldo kom ik tot de conclusie dat niet alleen Freudenthal maar nu ook David Tall een neiging hadden / hebben om Van Hiele in het hokje van de meetkunde te plaatsen, terwijl de Van Hieles juist in relatie tot Piaget een algemene theorie presenteerden met meetkunde slechts als voorbeeld.
>
> Ik houd me aanbevolen voor een reactie.
>
> Met vriendelijke groet,
>
> Thomas Cool / Thomas Colignatus
> Econometrist en leraar wiskunde
> Scheveningen

LB-VG's reply (Dutch original):

> From: Bastide-van Gemert, S la
> To: Thomas Cool / Thomas Colignatus

Subject: RE: N.a.v. uw proefschrift, Hoofdstuk 7, Van Hiele niveaux
Date: Tue, 19 Aug 2014

Geachte mijnheer Cool,

Dank voor uw e-mail. Ik heb nu helaas geen mogelijkheid inhoudelijk te reageren, maar wens u veel succes met uw artikel.

Met vriendelijke groet,

Sacha la Bastide

## *Appendix D: Applied mathematics in the autobiography*

Now that I am translating some quotes in 2015, I can supplement with what Freudenthal says in his autobiography about his view on the link between applied mathematics and "realistic mathematics eduation" (RME). In an autobiography an author might be excused for relating how things worked out for himself over time. Still:

- It would have been more efficient for him to refer to Van Hiele (1957), leaving out the suggestion as if he himself had a contribution that grew better over time.
- For an empirical researcher it would be natural to refer to studies that confirmed or rejected findings, instead of pontificating on personal views.

Since the Dutch text has three verbs jointly, I include an emphasis for reading:

> (a) "There is a field of 'applied mathematics'. I sensed at an early moment that teaching applied mathematics is not the way to get *learning to apply* mathematics, but only with the levels and the reflection it dawned on me why." [26]
>
> (b) "And thus it continues in traditional education in mathematics: *learning mathematics to apply it later* – a process that is didactically counterproductive. Most students aren't served by this even though the applications are the rationale of their learning mathematics – which they can't handle anyway. Here again we see the error of starting at a level to which you should grow instead." [27]
>
> (c) "How it must be done, I formulated in the following way, making it more precise over time: the reality in which you want to apply mathematics you must use first as source for the mathematics that you want to use in it. No application to it afterwards – yes, this too – but first of all investigate the field of application, even mathematise, nonconsciously, consciously and reflecting. In that manner the mathematics that you want to apply arises in the reality. This was the historical road, the road that you must also allow the learner to take, a stimulating allowance." [28]

Freudenthal *says* that he wants some distance to applied mathematics but *in practice* his method still relies on applied mathematics. He is familiar with the notion of application and not with the notion of teaching. As Van Hiele stated in the Alberts & Kaenders (2005) interview: Freudenthal may not really have understood the theory of levels. Not understanding something doesn't make the intellectual theft less so.

---

[26] Freudenthal (1987:357): "Er bestaat een vak 'toegepaste wiskunde'. Ik heb vroegtijdig aangevoeld dat het onderwijzen van toegepaste wiskunde niet de weg is om wiskunde te doen leren toepassen, maar pas met de niveaus en de reflectie werd me duidelijk waarom."

[27] "En zo gaat het door in het traditioneel wiskunde-onderwijs: wiskunde leren om achteraf toe te passen - een didactisch averechtse procedure, waarmee de meesten niet gediend zijn, ook al zijn die toepassingen - die ze toch niet aankunnen - dan het rationale van hun wiskunde-leren. Ook hier weer het euvel van het instappen op een niveau waar je eerst naar toe hoort te groeien."

[28] "Hoe het moet, heb ik zo geformuleerd en dan steeds scherper: de realiteit waar je wiskunde in wilt toepassen moet je allereerst als bron gebruiken voor die wiskunde die je erin wilt toepassen. Geen toepassen achteraf - ja dit ook - maar allereerst het gebied van toepassing verkennen, zelf mathematiseren, onbewust, bewust en reflecterend. Zo ontstaat in de realiteit de wiskunde die je daar wilt toepassen. Zo was historisch de gang van zaken, de weg die je ook de lerende moet toestaan om te bewandelen, stimulerend toestaan."

Ad (a): Freudenthal questions the effectiveness of standard teaching methods of applied mathematics *for learning mathematics*: but this is a category mistake to start with.

Courses in applied mathematics are *not intended* as courses to learn mathematics. The proper notion is that one must first understand mathematics before you can apply it. Throwing a ball is not applied mathematics but throwing a ball.

Like there are levels in mathematical insight, there are also levels in applied mathematics. You must still must understand mathematics at some level before you can apply it at that level. It is confusing to reject this (supposedly with the argument that competence at math level 1 does not allow application at level 2).

Ad (b): He creates a false opposition of RME to traditional education in mathematics (Euclid's *Elements*). Van Hiele already solved the confusion w.r.t. Euclid. The true opposition of RME is with the proposal by Van Hiele. If the views of VH and F would be the same then F can refer and that's it. By not referring he inserts the suggestion of a novel contribution by himself.

Ad (c): Apparently the frame of applied mathematics is so strong for Freudenthal – as can happen in academic Departments of Mathematics who only distinguish mathematics and applied mathematics (and who thus regard education as applied math too) – that he refers to application repeatedly. (We can also recognise his interest in history of mathematics, as a frame for teaching mathematics – nowadays seen as an element within mathematics education when it provides perspective.) Yet, when we acknowledge that Van Hiele had already discarded with the traditional ways (Euclid), it would have sufficed for Freudenthal to say that he followed Van Hiele in going from *concrete to abstract*, and prevent the confusion about application.

Given his repeated referral to application, his approach almost reads as a reform of education in applied mathematics. His text suggests that he wants to create some distance from applied mathematics, but by referring to it repeatedly he brings it closer. His repeated referral also is suggestive of some novelty w.r.t. Van Hiele.

Freudenthal likely would agree that the proper education of mathematics does not concern **a set of applications** but concerns **learning to think in general** – that would facilitate such applications too. Still, by deviating from the Van Hiele approach of *concrete to abstract*, he caused his followers to bring in all kinds of contexts within mathematics education that distracted from proper education of mathematics. (Kids learning arithmetic should not be distracted all the time by clocks and pizzas.)

## *Appendix E: Mathematical abolition of fractions (traditional notation)*

Consider the (crooked so-called "mathematical") division of mixed fractions: 2½ / 3⅓ .

Van der Plas (2008, 2009) reports that this kind of question hardly occurs in Dutch "realistic mathematics" textbooks. Traditional mathematicians like Hung-Hsi Wu in the USA want to see a lot of practice on this again. In handwriting by students: 2 ½ / 3 ⅓ = 1 / 3 ⅓ = 1 / 9. Van de Craats in Holland is aware of the issue and suggests using 5 / 2 but then loses the location on the number line.

The mathematical meaning of the inverse $x$^$H$ or $x^H$ is that: $x\, x^H = 1$ (for $x \neq 0$).

On the calculator we find a numerical approximation of $x^H$ by: $(x)$^$(-1)$.

A relation for exponents is: $x = (x^H)^H$.

Above expression becomes in neoclassical mathematics: $(2 + 2^H)(3 + 3^H)^H$

For pupils this is a new notation. Who is used to it may take longer strides. For now, we take small steps. A classical operation may be at least as long. Properties of *H* ("eta") are stated in the book "A child wants nice and no mean numbers" (2015) with some observations that might be useful for elementary school (I have no degree on that area). Eventually, students must learn to handle exponents. The following might seem complex but eventually it will be faster and more insightful. Obviously, this is only an expectation, and it must be checked with pupils whether this expectation is corroborated. It are the pupils who determine what works.

For ease of comparison, the following two schemes are put onto a single page.

Operations by steps:

| | |
|---|---|
| $(2 + 2^H) \; (3 + 3^H)^H$ | Objective: simplify as much as possible |
| $(2 \; (2 \, 2^H) + 2^H) \; (3 \, (3 \, 3^H) + 3^H)^H$ | Use $x \, x^H = 1$ |
| $(4 \, 2^H + 2^H) \; (9 \, 3^H + 3^H)^H$ | Multiply |
| $(4 + 1) \, 2^H \; ((9 + 1) \, 3^H)^H$ | Take factors out of the brackets |
| $5 \, 2^H \; 10^H \; 3$ | Add up weights |
| $5 \, 2^H \; (5 \, 2)^H \; 3$ | Factorise |
| $5 \, 5^H \; 2^H \; 2^H \; 3$ | Use $x \, x^H = 1$ |
| $3 \, 4^H$ | 3 per 4 |

For comparison: the crooked manner of traditional so-called "mathematics", in which the notation determines what must be done, and in which you don't just denote what you are doing.

| | |
|---|---|
| 2½ / 3⅓ | Dangerous notation |
| (2 + ½) / (3 + ⅓) | Turn it into mathematics, otherwise it doesn't work |
| (2 (2 ½) + ½) / (3 (3 ⅓) + ⅓) | Use $x / x = 1$ (notation causes a problem) |
| (4 ½ + ½) / (9 ⅓ + ⅓) | Multiply |
| (4 / 2 + ½) / (9 / 3 + ⅓) | Same denominators (extra concept) |
| (5 / 2) / (10 / 3) | Add up weights (5 / 2: number or operation ?) |
| (5 / 2) (3 / 10) | Division by fraction is multiplication by inverted |
| (5 3) / (2 10) | Multiply |
| 3 / (2 2) | Factorise and eliminate equal terms |
| 3 / 4 | three-fourths (abuse of rank order number "fourth") |